\documentclass[a4paper,twoside]{article}

\usepackage{reconstructible}

\author{Mireille Boutin and Gregor Kemper
}

\title{On Reconstructing $n$-Point Configurations from the
  Distribution of Distances or Areas}

\begin{document}

\maketitle

\begin{abstract}
  One way to characterize configurations of points up to congruence is
  by considering the distribution of all mutual distances between
  points.  This paper deals with the question if point configurations
  are uniquely determined by this distribution. After giving some
  counterexamples, we prove that this is the case for the vast
  majority of configurations.
  
  In the second part of the paper, the distribution of areas of
  sub-triangles is used for characterizing point configurations.  Again
  it turns out that most configurations are reconstructible from the
  distribution of areas, though there are counterexamples.
\end{abstract}

%\tableofcontents

\section*{Introduction} \label{0sIntro}

In this paper, 
we study a type of shape representation
which attempts
to combine both the approaches of invariant theory and statistics. 
We consider the problem of characterizing the {\em shape} or, 
more generally,
the {\em geometry} of a configuration of points.
More precisely, we are interested in finding a good representation for 
configurations of points in a vector space
modulo the action of a Lie group $G$.
The solution we investigate 
consists in using distributions of invariants of the action of $G$.

Our main motivation comes from applications in computer vision.
A central problem in image understanding
is that of identifying 
objects from a picture.
In that problem, 
one must take into account that
variations in the position of the object or 
in the parameters of the camera 
induce variations in the image 
which correspond to group transformations 
that need to be moded out 
in order to establish the correspondence between two 
pictures of the same object.

The obvious way to obtain 
image features which are not affected by the action of the group
is to use invariants of the group action.
However, in order to be able to positively identify any object,
we need to find a set of invariants 
whose values {\em completely} characterize 
the image of the object up to the action of the group.
In other words, we need to find a set of invariants such that
two images are in the same orbit
{\em if and only if }
the values of these invariants evaluated on the two images
are the same.
Such invariants are called {\em separating}
because they can be used to separate  the orbits. 
In traditional approaches to object recognition
(see for example \mycite{Mundy:Zisserman}), 
this method is commonly used.

In the following, we address the case of shapes defined by a finite
set of points.  This is actually an important case for applications.
Indeed for many reasons (e.g.~the amount of noise or the nature of the
data) it is common to represent an object of interest by a finite set
of points called {\em landmarks}.  For example, landmarks can be
defined by salient features on the boundary of the
image of the object. Specifically, one might think of minutiae in
fingerprints, corners on edges of archaeological sherds, or stellar
constellations.  In order to recognize the object, one thus needs to
characterize the point configuration given by the landmarks up to the
action of the group.

Given a Lie group $G$ acting on a vector space $V$
and two sets of n points
$P_1,\ldots,P_n$ and $\bar{P}_1,\ldots,\bar{P}_n\in V$,
we want to be able to determine whether there exists
$g\in G$ and a permutation $\pi \in S_n$ 
(since, a priori, 
we don't know whether the points are labeled in correspondence) 
such that
\[
g ( P_i )=\bar{P}_{\pi(i)}, \text{ for all }i=1,\ldots,n.
\]
In applications, we are often interested in pictures, so
$V$ is usually ${\mathbb R}^2$ or ${\mathbb R}^3$
and
the Lie group $G$ is typically 
a subgroup of the projective group 
and depends on how the picture of the object was taken.
Examples of important groups include
$\AO(2)$, the group of rigid motions in the plane 
(rotations, reflections and translations, sometime also denoted by $E(2)$),
and $A(2)$, the group of affine transformations in the plane, i.e.~all translations and linear maps with determinant $\pm 1$.

In principle, this problem can indeed be solved 
using invariants.
If we assume that the points are distinguishable
so we know how to correctly label them, 
then all we need to do is to find a set of separating invariants 
of the diagonal action of $G$ on  $V^n$,
\[
g\cdot (Q_1,\ldots,Q_n) =(g ( Q_1),\ldots, g (Q_n) )\text{ for all }g\in G,\text{ and all } Q_1,\ldots,Q_n\in V .
\]
For example, 
if $G=\AO(2)$ the group of Euclidean transformations in the plane 
then two sets of landmarks 
$P_1,\ldots,P_n$ and $\bar{P}_1,\ldots,\bar{P}_n$ 
(labeled in correspondence)
belong 
to the same orbit under the action of $\AO(2)$
if and only if 
all their pairwise distances $d(P_i,P_j)=d(\bar{P}_i,\bar{P}_j)$ 
are the same for all $i,j=1,\ldots, n$.
So the shape of the set of labeled landmarks $P_1,\ldots ,P_n$
is completely characterized 
by the value of the pairwise (labeled) distances between the landmarks.

However, in  most applications
the point correspondence is unknown
so things are more complicated,
especially when the number of points $n$ is big.
Indeed, labeling the points is a non-trivial task
which, although feasible, takes time. 
(See for example \mycite{HartleyZisserman} 
for an easy exposition of some existing methods.)
And the bigger the number of points, the longer it takes.
We would thus prefer to simply skip the labeling step.
So, can we, instead, find separating invariants 
of the action of $\AO(2)\times S_n$?

The answer to this question is, of course, yes.
For example, in the case $n=3$,
instead of distances one can use
the following symmetric functions of the distances,
\begin{eqnarray*}
f_1(P_1,P_2,P_3) &=& d(P_1,P_2)+d(P_1,P_3)+d(P_2,P_3),\\
f_2(P_1,P_2,P_3) &=& d(P_1,P_2)d(P_2,P_3)+d(P_1,P_2)d(P_1,P_3)+d(P_1,P_3)d(P_2,P_3),\\
f_3(P_1,P_2,P_3)&=& d(P_1,P_2)d(P_1,P_3)d(P_2,P_3).
\end{eqnarray*}
These are separating invariants of the action of $\AO(2)\times S_3$ on
$\left( {\mathbb R}^2\right)^3$. Continuing in this way, we can try to
find expressions in the distances $d(P_1,P_2)$, $d(P_1,P_3)$,
$d(P_1,P_4)$, $d(P_2,P_3)$ $d(P_2,P_4)$, and $d(P_3,P_4)$, which are
invariant under the action of $S_4$ by permuting the $P_i$, and which
form a generating (or at least separating) subset of all such
invariants.  But notice that the elementary symmetric functions in the
distances will not qualify anymore, since these are the invariants
under the action of $S_6$ instead of $S_4$. Thus this approach
requires a fresh computation of invariants for each value of~$n$.

The $S_n$-invariants needed here are often called {\em graph
  invariants}, and have been studied in a graph theoretical context by
various authors, e.g.  \mycite{Nicolas.Thiery:c}, \mycite{Pouzet}, and
\mycite{ACG}.  \mycite{ACG} calculated a generating set of graph
invariants for $n = 4$, obtaining a minimal set of~9 invariants. But
for $n = 5$ the computation of graph invariants is already very hard
and stood as a challenge problem for a while (see~[\citenumber{ACG},
\citenumber{Nicolas.Thiery:c}]) until the computation was done by the
second author (see \mycite[p.~221]{Derksen:Kemper}). The minimal
generating set for $n = 5$ contains~56 invariants, and storing them
takes several MBytes of memory. For $n \ge 6$ the computation is
presently not feasible. This clearly shows that the approach of using
graph invariants is far from practical. Apart from their number and
the difficulties of computing them, they cannot be used in practice
for questions of robustness, since high degree polynomials vary
immensely when small variations in the points $P_1,\ldots,P_n$ are
introduced.  We thus need to find better invariants than graph
invariants; we need invariants that not only separate the orbits of
the action of $G\times S_n$ but that are also robust and simple to
compute.

We were inspired by looking at what engineers 
do in practice.
In order to identify images of the same object, 
they often drop the separation requirement
and simply look for invariant features of the image
of which they compare the distribution.
The distribution of the pairwise distances of a set of points
is obviously invariant under a relabeling of the point. 
It is also much more robust than 
a set of polynomial functions of the pairwise distances.
In addition, it is not too complicated to compute 
and very easy to manipulate.

So we asked ourselves if the distribution of distances of a set of
points is actually also a separating invariant and thus completely
characterizes point configurations up to rigid motions. In other
words, can an $n$-point configuration be reconstructed uniquely (up to
the labeling of the points and up to rigid motions) from the
distribution of distances?  It turns out that this is {\em false} in
general, as we demonstrate with counterexamples.  But fortunately,
counterexamples are rare, in a sense to be explained shortly.  This is
the contents of our first main result (\tref{1tReconstructible}).
Moreover, it is true locally, i.e.~the shape of $n$-point
configurations that are close enough can be compared using their
distribution of invariants. We also explore methods to verify
reconstructibility for particular configurations. Most of the results
for the case of distances in the real plane naturally extend to any
vector space with a non-degenerate quadratic form over a field of
characteristic not equal to 2.  We shall thus simply treat this
general case in the first part of this paper.

In the second part, we attempt to characterize point configurations up
to the action of the equi-affine group $A(2)$ and, again, the
symmetric group $S_n$. This action is relevant in computer vision
since, up to a scale factor, it adequately approximates what happens
to the camera image of a very distant planar object as it is rotated
and translated in three-dimensional space.  As above, there are
obvious invariants for separating orbits under $A(2)$. These are the
areas of triangles spanned by a selection of three of the~$n$ points.
As before, we attempt to separate $S_n$-orbits by considering the {\em
  distribution} of all these areas. We obtain results which are
completely analogous to those in the first section: There are examples
of configurations which cannot be reconstructed (up to the action of
$A(2) \times S_n$) from the distribution of areas; but a dense open
subset of configurations are reconstructible in this sense (see
\tref{2tReconstructible}). We believe that for most purposes in
computer vision, this is a satisfactory result. Again our results
generalize to configurations in any dimension and to any ground field.

Let us emphasize here that the use of computer algebra systems played
a vital role in the preparation of this paper. In particular,
Magma~[\citenumber{magma}] was an indispensable tool. For example, the
first example of an $n$-point configuration which is not
reconstructible from distances was the upshot of a prolonged Magma
session. The examples in Sections~\ref{2sCounter} and~\ref{2sCombine}
were constructed with the help of Magma and Maple~[\citenumber{map}].
But also the proof of \tref{1tReconstructible} was inspired by sample
computations in Magma.

\section{Reconstruction from distances} \label{1sDistances}

An $n$-point configuration is a tuple of points $P_1 \upto P_n \in
\RR^m$.  To an $n$-point configuration we associate the squared
(Euclidean) distances $d_{i,j}$ between each pair of points $P_i$ and
$P_j$, and then consider the {\em distribution} of distances, i.e.~the
relative frequencies of the value of the distances.  In other words,
the distribution of distances of an $n$-point configuration tells us
how many times each distance occurs relative to the total number of
distances.  This means that, for $n$ fixed, the distribution of
distances is given by the set of the numbers $d_{i,j}$ possibly with
multiplicities if some distances occur several times.  So considering
the distribution of distances of an $n$-point configuration is
equivalent to considering the polynomial
\[
F_{P_1 \upto P_n}(X) := \prod_{1 \leq i < j \leq n} (X - d_{i,j}).
\]
In order to better visualize the information contained 
in a distribution of distances, 
one can plot a histogram of the distances, 
i.e.~one can group the data into bins of a fixed size 
and count how many distances lie in each bin. 
Figure~\ref{fig:rectangle1}, \ref{fig:annulus1} and \ref{fig:tworectangles1}
show examples of $n$-point configurations in the plane
together with a histogram of their distances.

\Figure{  
\begin{tabular}{cc}   
\includegraphics[height=3.0cm]{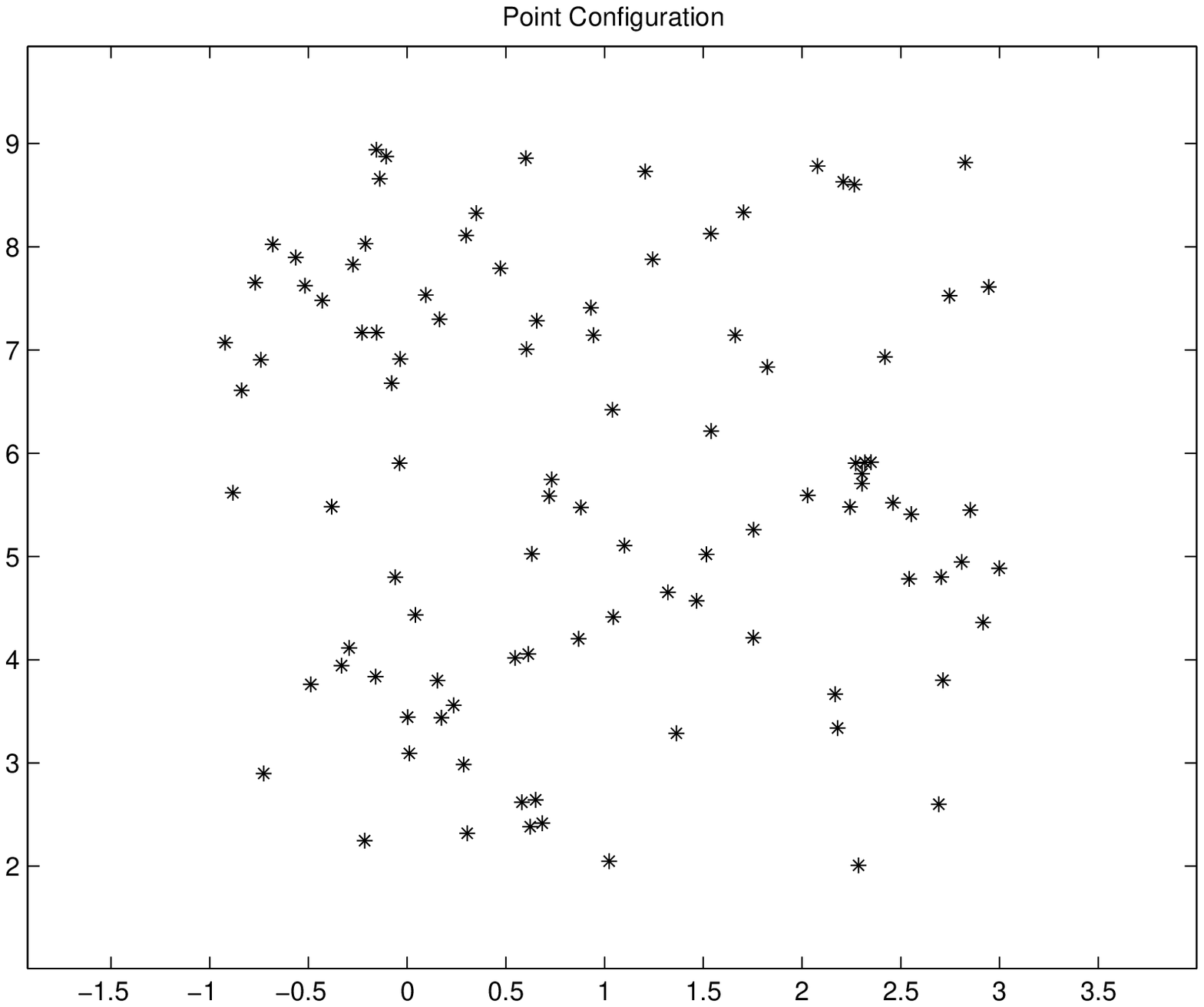} &  
\includegraphics[height=3.0cm]{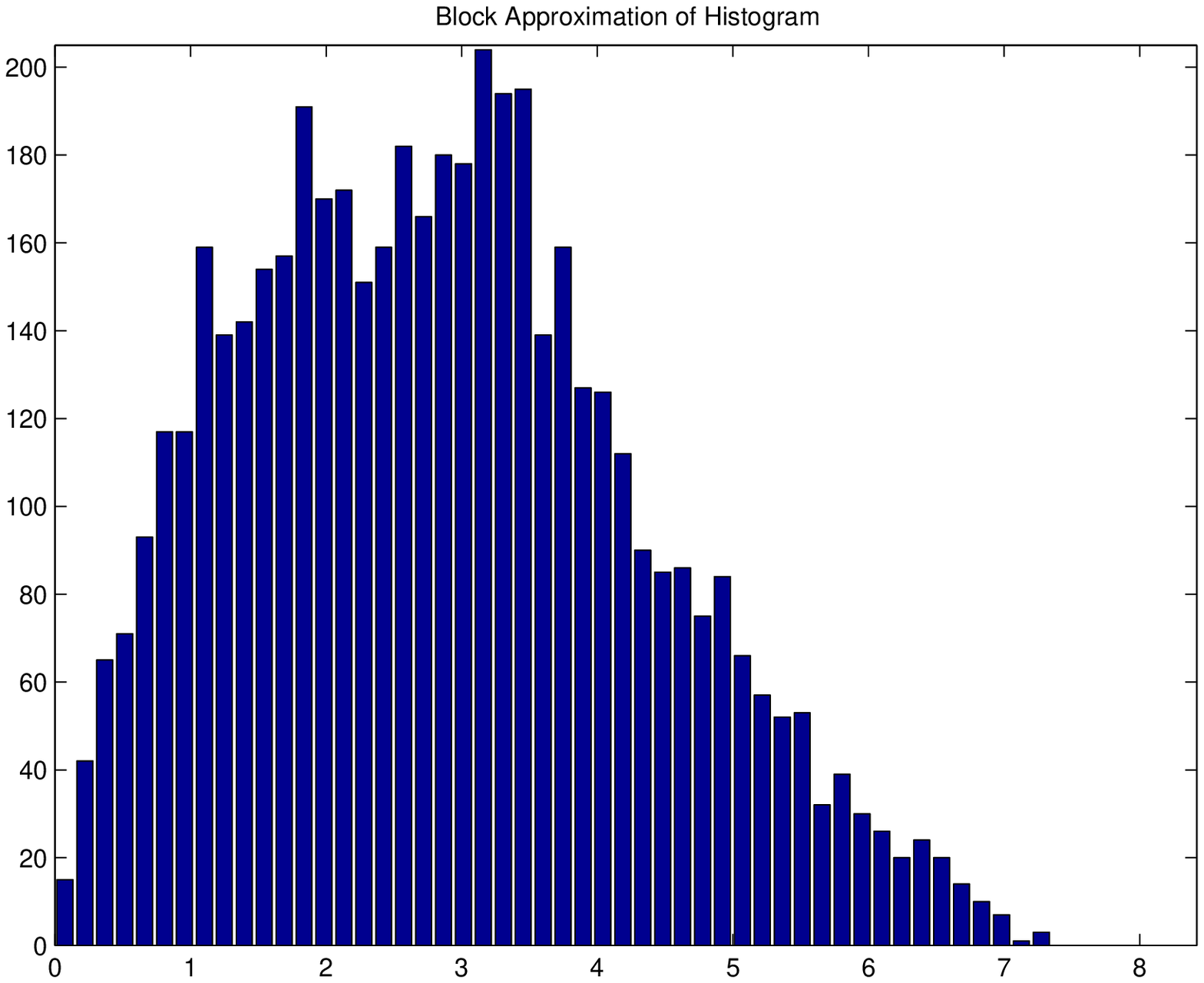} \\
a & b  \\    
\end{tabular}    
}{a) A $100$-point configuration, b) Histogram of distances with bin size 0.1470}{fig:rectangle1} 

\Figure{  
\begin{tabular}{cc}   
\includegraphics[height=3.0cm]{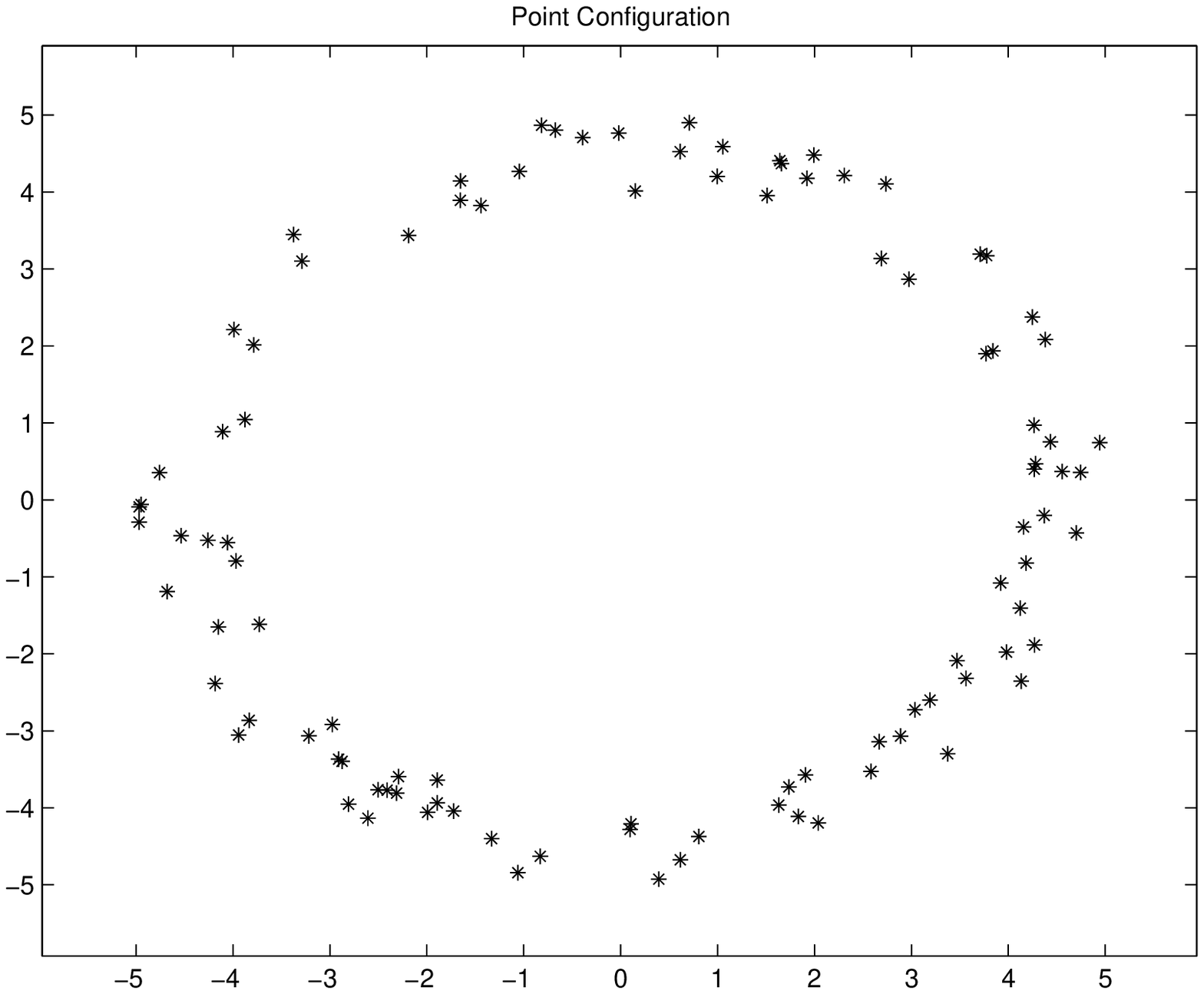} &  
\includegraphics[height=3.0cm]{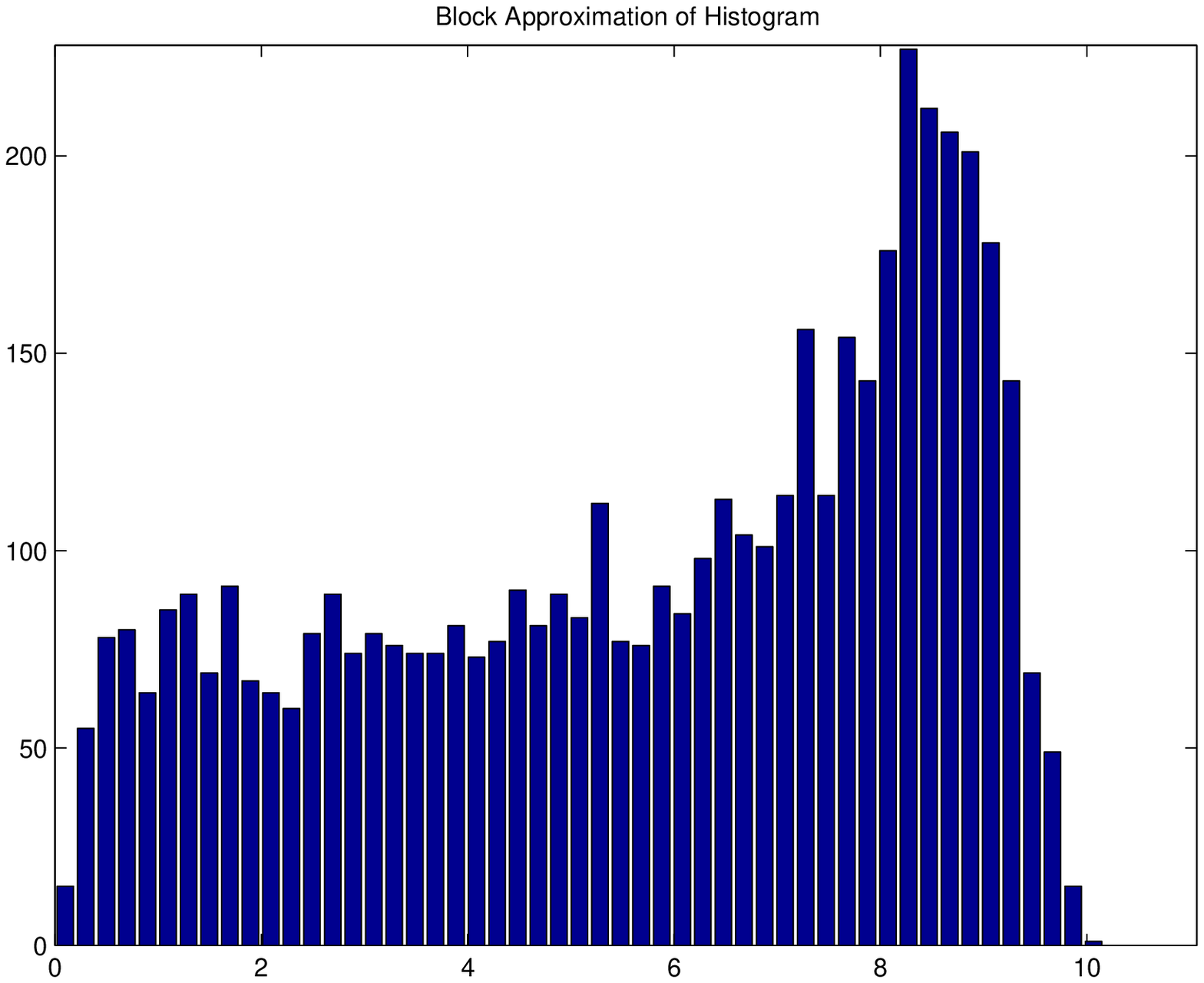} \\
a & b  \\    
\end{tabular}    
}{a) A $100$-point configuration, b) Histogram of distances with bin size 0.1993}{fig:annulus1} 

\Figure{  
\begin{tabular}{cc}   
\includegraphics[height=3.0cm]{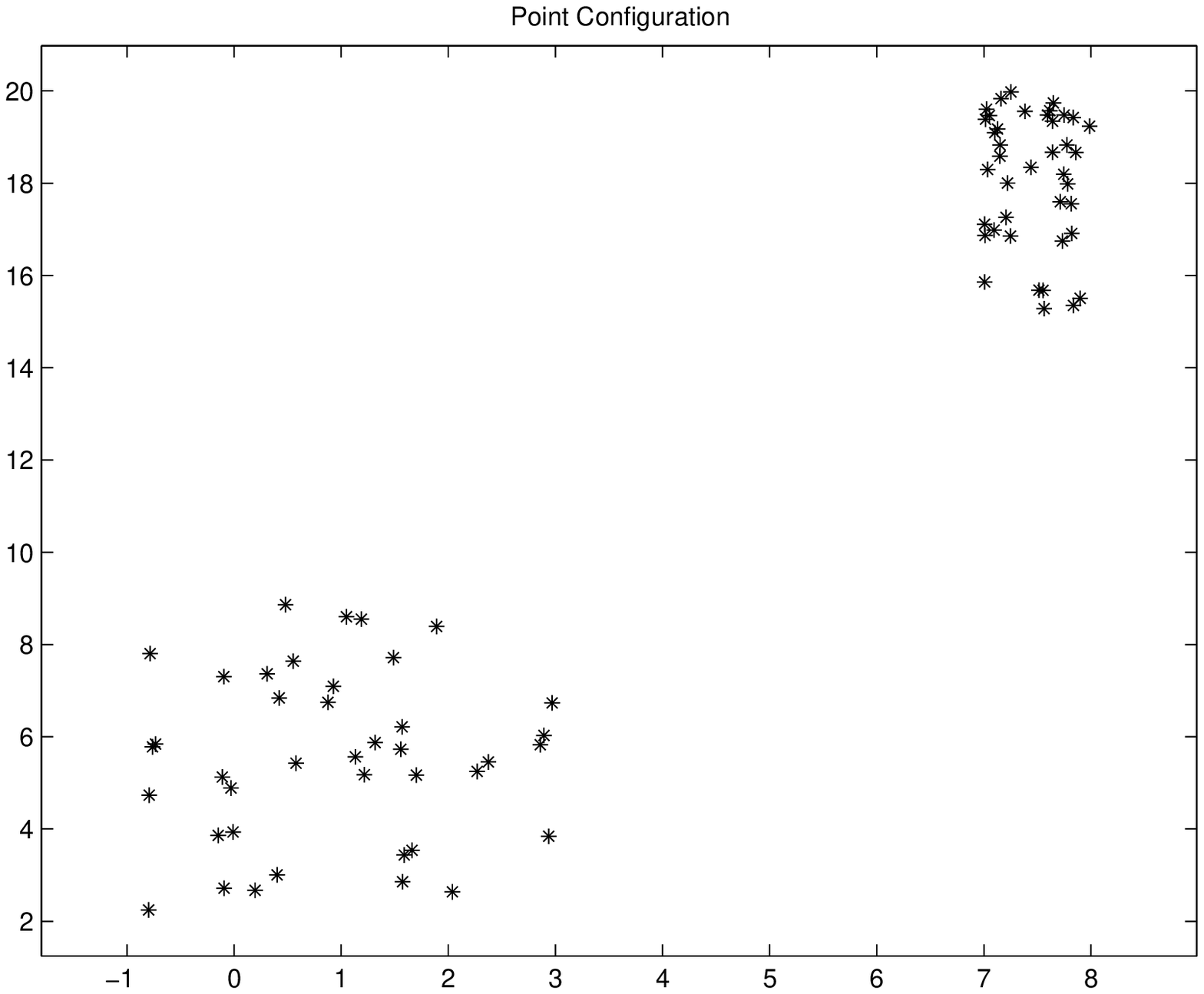} &  
\includegraphics[height=3.0cm]{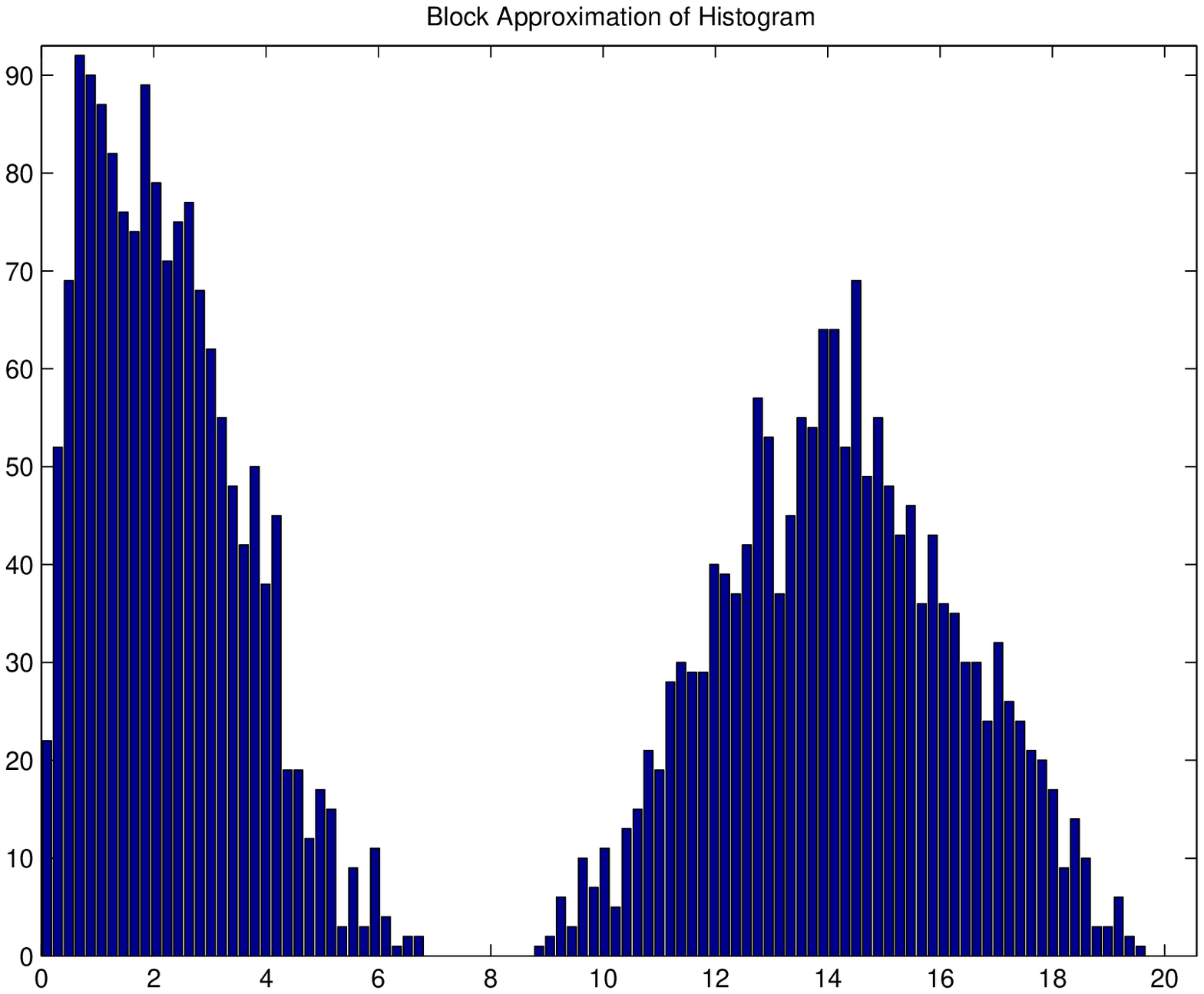} \\
a & b  \\    
\end{tabular}    
}{a) An $80$-point configuration, b) Histogram of distances with bin size 0.1947}{fig:tworectangles1}

Clearly the distribution of distances is {\em invariant} under
permutations of the points and under the (simultaneous) action of the
Euclidean group. The question is whether an $n$-point configuration
can be reconstructed from its distribution of distances.

\begin{defi} \label{1dReconstructible}
  An $n$-point configuration $P_1 \upto P_n \in \RR^m$ is called
  \df{reconstructible from distances} if the following holds. If $Q_1
  \upto Q_n$ is another $n$-point configuration with $F_{P_1 \upto
  P_n}(X) = F_{Q_1 \upto Q_n}(X)$, then there exists a permutation
  $\pi \in S_n$ and a Euclidean transformation $g \in \AO_m(\RR)$ such
  that $g(P_{\pi(i)}) = Q_i$ holds for all~$i$.
\end{defi}

The notion of reconstructibility from distances generalizes naturally
to any vector space with a non-degenerate quadratic form over a field
of characteristic not equal to~2.

\subsection{Non-reconstructible configurations} \label{1sCounter}

It is clear that in two-dimensional Euclidean space all triangles are
reconstructible from distances, and the same is true for all 2-point
configurations. So the quest for examples of non-reconstructible
$n$-point configurations becomes interesting for $n \ge 4$.
Figure~\ref{1fCounter} shows such an example. We have put the
(non-squared) distances along the lines connecting pairs of points.
Note that the upper point in the first configuration is moved
diagonally downward to obtain the second configuration, while the
other three points remain inert.

\newcommand{\point}{\circle*{2.5}}
\newcommand{\attach}[1]{\scriptsize $#1$}
\Figure{
  \fbox{
    \begin{picture}(130,90)
      \put(5,45){\point}
      \put(95,75){\point}
      \put(95,15){\point}
      \put(125,45){\point}
      \put(5,45){\line(1,0){120}}
      \put(42,47){\attach{4}}
      \put(5,45){\line(3,1){90}}
      \put(61,69){\attach{\sqrt{10}}}
      \put(5,45){\line(3,-1){90}}
      \put(43,33){\attach{\sqrt{10}}}
      \put(125,45){\line(-1,1){30}}
      \put(109,61){\attach{\sqrt{2}}}
      \put(125,45){\line(-1,-1){30}}
      \put(101,32){\attach{\sqrt{2}}}
      \put(95,75){\line(0,-1){60}}
      \put(90,59){\attach{2}}
      \dottedline[\circle*{0.4}]{2}(107,3)(23,87)
      \put(38.15,71.85){\circle*{1}}
    \end{picture}
  }
  \hspace{11mm}
  \fbox{
    \begin{picture}(130,90)
      \put(5,45){\point}
      \put(35,15){\point}
      \put(95,15){\point}
      \put(125,45){\point}
      \put(5,45){\line(1,0){120}}
      \put(43,47){\attach{4}}
      \put(5,45){\line(3,-1){90}}
      \put(41,34){\attach{\sqrt{10}}}
      \put(5,45){\line(1,-1){30}}
      \put(8,23){\attach{\sqrt{2}}}
      \put(125,45){\line(-3,-1){90}}
      \put(83,36){\attach{\sqrt{10}}}
      \put(125,45){\line(-1,-1){30}}
      \put(106,23){\attach{\sqrt{2}}}
      \put(35,15){\line(1,0){60}}
      \put(63,8){\attach{2}}
      \dottedline[\circle*{0.4}]{2}(107,3)(23,87)
      \put(38.15,71.85){\circle*{1}}
    \end{picture}
  }
}
{Two 4-point configurations with the same distribution of distances}
{1fCounter}

Further examples can be constructed by adding an arbitrary number of
additional points on the dotted line and at the same position in both
configurations (such as the slightly thicker dot in each picture).
Thus we get examples of non-reconstructible $n$-point configurations
for any $n \ge 4$. By embedding these into a space of higher
dimension, we also get examples in any dimension $m \ge 2$. The fact
that we can add points at arbitrary positions on the dotted line shows
that the symmetry of the configuration is not responsible for the fact
that it is not reconstructible.

\subsection{Relation-preserving permutations} \label{1sRelation}

Let $K$ be a field of characteristic not equal to~2 ($K = \CC$ and $K
= \RR$ will be the most important examples). Let $V$ be an
$m$-dimensional vector space over $K$ with a non-degenerate symmetric
bilinear form $\langle \cdot,\cdot \rangle$. With a suitable choice of
a basis, this form is given by $\langle (x_1 \upto x_m),(y_1 \upto
y_m) \rangle = \sum_{k=1}^m a_k x_k y_k$ with $a_k \in K \setminus
\{0\}$. If $v_1 \upto v_n$ are vectors in $V$, 
then the Gram matrix \linebreak
$\left(\langle v_i,v_j\rangle\right)_{i,j = 1 \upto n}$ has rank at
most~$m$, hence the $(m+1) \times (m+1)$-minors are zero. By the
following well-known proposition, this gives all relations between the
scalar products of~$n$ vectors. Part~(b) gives the relations between
the distances between~$n$ points. In fact, \pref{1pRelations}(a) is
the ``second fundamental theorem'' of invariant theory of orthogonal
groups.

\begin{prop} \label{1pRelations}
  Let $x_{i,k}$ be indeterminates over $K$ ($i = 1 \upto n$, $k = 1
  \upto m$).
  \begin{enumerate}
  \item Let $s_{i,j}$ be further indeterminates ($1 \leq i \leq j \leq
    n$). Then the kernel of the map
    \[
    K[s_{1,1} \upto s_{n,n}] \to K[x_{1,1} \upto x_{n,m}],\ s_{i,j}
    \mapsto \sum_{k=1}^m a_k x_{i,k} x_{j,k}
    \]
    is generated (as an ideal) by the $(m+1) \times (m+1)$-minors of
    the matrix $\left(s_{i,j}\right)_{i,j = 1 \upto n}$, where we set
    $s_{i,j} := s_{j,i}$ for $i > j$.
  \item Let $D_{i,j}$ be indeterminates ($1 \leq i < j \leq n$). Then
    the kernel of the map
    \[
    K[D_{1,2} \upto D_{n-1,n}] \to K[x_{1,1} \upto x_{n,m}],\ D_{i,j}
    \mapsto \sum_{k=1}^m a_k \left(x_{i,k} - x_{j,k}\right)^2
    \]
    is generated (as an ideal) by the $(m+1) \times (m+1)$-minors of
    the matrix
    \begin{equation} \label{1eqD}
      {\cal D} = \left(D_{i,j} - D_{i,n} - D_{j,n}\right)_{i,j = 1
        \upto n-1},
    \end{equation}
    where we set $D_{i,i} := 0$ and $D_{i,j} := D_{j,i}$ for $i > j$.
  \end{enumerate}
\end{prop}

\begin{proof}
  For part~(a), see \mycite{Weyl} or
  \mycite[Theorem~5.7]{deConcini:Procesi} (the latter reference takes
  care of the positive characteristic case). Part~(b) follows from~(a) 
  since for points $P_1 \upto P_n \in V$ we have
  \begin{equation} \label{1eqSD}
    \langle P_i - P_n,P_j - P_n\rangle = \frac{1}{2} \left(\langle P_i
      - P_n,P_i - P_n \rangle + \langle P_j - P_n,P_j - P_n \rangle -
      \langle P_i - P_j,P_i - P_j \rangle\strut\right).
  \end{equation}
\end{proof}

We will now study monomials occurring in elements of the ideal given
by \pref{1pRelations}(b). From now on it is useful to use sets
$\{i,j\}$ as indices of the $d$'s rather than pairs $i,j$.

\begin{lemma} \label{1lMonomials}
  Let $K$ be a field of characteristic not equal to~2 and let
  $D_{\{i,j\}}$ be indeterminates ($i,j = 1 \upto n$, $i \ne j$). For
  an integer $r$ with $1 \leq r \leq n-1$ consider the ideal $I$
  generated by all $(r \times r)$-minors of the matrix ${\cal D} :=
  \left(D_{\{i,j\}} - D_{\{i,n\}} - D_{\{j,n\}}\right)_{i,j = 1 \upto
    n-1}$, where we set $D_{\{i,i\}} := 0$. Let $t = \prod_{\nu = 1}^r
  D_{\{i_\nu,j_\nu\}}$ be a monomial of degree~$r$. Then the following
  are equivalent:
  \begin{enumerate}
  \item The monomial~$t$ occurs in a polynomial from $I$.
  \item Every index from $\{1 \upto n\}$ occurs at most twice among
    the $i_\nu$ and $j_\nu$. More formally, for every $k \in \{1 \upto
    n\}$ we have $|\{\nu \mid i_\nu = k\}| + |\{\nu \mid j_\nu = k\}|
    \leq 2$.
  \end{enumerate}
\end{lemma}

\begin{proof}
  It follows from \pref{1pRelations}(b) that the ideal $I$ is stable
  under the natural action by the symmetric group $S_n$. Thus~$t$
  occurs in a polynomial from $I$ if and only if all images of~$t$
  occur.
  
  First assume that there exists a $k \in \{1 \upto n\}$ which occurs
  more than twice among the $i_\nu$ and $j_\nu$. By the previous
  remark we may assume $k = 1$. If $t$ occurs in a polynomial of $I$
  it must also occur in an $(r \times r)$-minor of $\cal D$ (since
  $\deg(t) = r$). But in order to obtain~$t$ as a monomial in an $(r
  \times r)$-minor, one has to choose the first row or the first
  column of $\cal D$ at least twice, since entries involving the
  index~1 only occur in the first row and column. But that is
  impossible. This proves that~(a) implies~(b).
  
  Now assume that~(b) is satisfied. Consider the graph $\cal G$ with
  vertices indexed $1 \upto r$, where the number of edges between
  vertex~$\nu$ and~$\mu$ is $|\{i_\nu,j_\nu\} \cap \{i_\mu,j_\mu\}|$,
  i.e., the number of indices shared by the $\nu$-th and $\mu$-th
  indeterminate in~$t$. By the hypothesis~(b) every vertex is
  connected to at most two others, hence every connected component of
  $\cal G$ is a line (including the case of an unconnected vertex) or
  a loop (including the case of a loop of two vertices corresponding
  to indeterminates $D_{\{i_\nu,j_\nu\}}$ and $D_{\{i_\mu,j_\mu\}}$
  which are equal). By renumbering, we may assume that the first
  connected component is given by the first~$m$ vertices. By the
  remark at the beginning of the proof, we may further assume that the
  first~$m$ indeterminates in~$t$ are $D_{\{1,2\}},D_{\{2,3\}} \upto
  D_{\{m,m+1\}}$ (forming a line in $\cal G$) or
  $D_{\{1,2\}},D_{\{2,3\}} \upto D_{\{m-1,m\}},D_{\{1,m\}}$ (a loop).
  Since $m \leq r \leq n-1$, it can only happen in the first case that
  the index~$n$ is involved in these indeterminates. Thus if~$n$ is
  involved, then $m = n - 1$ and $t = \prod_{\nu=1}^{n-1}
  D_{\{\nu,\nu+1\}}$. It is easily seen that in this case~$t$ occurs
  in $\det({\cal D})$ with coefficient $2 \cdot (-1)^{n-1}$. Having
  settled this case, we may assume that $m < n-1$. We proceed by
  induction on the number of connected components of $\cal G$.
  
  First assume that the first component is a loop. We wish to build an
  $(r \times r)$-submatrix of $\cal D$ whose determinant contains~$t$
  as a monomial. To this end, we start by choosing the first~$m$ rows
  and the first~$m$ columns from $\cal D$. Temporarily setting all
  $D_{\{i,n\}} := 0$, we obtain a matrix ${\cal D}'$ with
  \[
  {\cal D}'|_{_{D_{\{i,n\}} = 0}} =
  \begin{pmatrix}
    0 & D_{\{1,2\}} & \cdots & D_{\{1,m-1\}} & D_{\{1,m\}} \\
    D_{\{1,2\}} & 0 & \cdots & D_{\{2,m-1\}} & D_{\{2,m\}} \\
    \vdots & & \ddots & & \vdots \\
    D_{\{1,m-1\}} & D_{\{2,m-1\}} & \cdots & 0 & D_{\{m-1,m\}} \\
    D_{\{1,m\}} & D_{\{2,m\}} & \cdots & D_{\{m-1,m\}} & 0
  \end{pmatrix}.
  \]
  Clearly the product $t' := D_{\{1,2\}} D_{\{2,3\}} \cdots
  D_{\{m-1,m\}} D_{\{1,m\}}$ occurs with coefficient $2 \cdot
  (-1)^{m-1}$ (or $-1$ if $m = 2$) in $\det({\cal D}')$. Since the
  first~$m$ indeterminates in~$t$ form a connected component in $\cal
  G$, the indeterminates in $t'' := t/t'$ involve none of the indices
  $1 \upto m$. Thus by induction we can choose $r-m$ rows, all below
  the $m$-th row, and $r-m$ columns, all right of the $m$-th column,
  such that $t''$ occurs as a monomial of the determinant of the
  corresponding submatrix ${\cal D}''$.  Finally, in order to get all
  of~$t = t' \cdot t''$ as a monomial in a minor, choose the rows and
  columns as in ${\cal D}''$ together with the first $m$ rows and
  columns. This yields a submatrix of $\cal D$ of block structure
  \[
  \begin{pmatrix}
    {\cal D}' & * \\
    * & {\cal D}''
  \end{pmatrix},
  \]
  where indeterminates $D_{\{i,j\}}$ with both indices $ \leq m$ only
  occur in ${\cal D}'$. Now clearly~$t$ occurs with non-zero
  coefficient in the determinant of this matrix.
  
  Let us treat the second case, so assume that the first component of
  $\cal G$ is a line $D_{\{1,2\}},D_{\{2,3\}} \linebreak\upto D_{\{m,m+1\}}$.
  Taking rows $1 \upto m$ and columns $2 \upto m+1$ yields a matrix
  ${\cal D}'$ with
  \[
  {\cal D}'|_{_{D_{\{i,n\}} = 0}} =
  \begin{pmatrix}
    D_{\{1,2\}} & D_{\{1,3\}} & \cdots & D_{\{1,m\}} & D_{\{1,m+1\}} \\
    0 & D_{\{2,3\}} & \cdots & D_{\{2,m\}} & D_{\{2,m+1\}} \\
    D_{\{2,3\}} & 0 & \cdots & D_{\{3,m\}} & D_{\{3,m+1\}} \\
    \vdots & & \ddots & & \vdots \\
    D_{\{2,m\}} & D_{\{3,m\}} & \cdots & 0 & D_{\{m,m+1\}} \\
  \end{pmatrix}.
  \]
  The product $t' := D_{\{1,2\}} D_{\{2,3\}} \cdots D_{\{m,m+1\}}$
  occurs with coefficient~1 in $\det({\cal D}')$. As above, the
  monomials in the remaining part $t'' := t/t'$ of~$t$ only involve
  indices strictly bigger than $m + 1$. Thus we may choose $r - m$
  rows and columns which are all below and right of the $(m+1)$-st,
  respectively, to form a submatrix ${\cal D}''$ which has $t''$ in
  its determinant. Again, putting together the rows and columns that
  we chose yields a submatrix with block structure as above. We see
  that also in this case~$t$ occurs as a monomial in an $(r \times
  r)$-minor of $\cal D$.
\end{proof}

If two $n$-point configurations have the same distribution of
distances, this means that the distances of both configurations
coincide up to some permutation.  But the permuted distances must
again satisfy the relations given by the ideal from
\pref{1pRelations}. Therefore it is crucial to determine how this
ideal behaves under permutations of the $D_{\{i,j\}}$. We show that
all permutations which preserve this ideal are in fact induced from
permutations of the~$n$ points. This provides the core of our
argument.

\begin{lemma} \label{1lPermutations}
  Let $K$ be a field of characteristic not equal to~2 and let
  $D_{\{i,j\}}$ be indeterminates ($i,j = 1 \upto n$, $i \ne j$). For
  an integer $r$ with $3 \leq r \leq n-1$ consider the ideal $I$
  generated by all $(r \times r)$-minors of the matrix ${\cal D} :=
  \left(D_{\{i,j\}} - D_{\{i,n\}} - D_{\{j,n\}}\right)_{i,j = 1 \upto
    n-1}$, where we set $D_{\{i,i\}} := 0$. Let~$\phi$ be a
  permutation of the $D_{\{i,j\}}$ which maps $I$ to itself. Then
  there exists a permutation $\pi \in S_n$ such that
  \[
  \phi(D_{\{i,j\}}) = D_{\{\pi(i),\pi(j)\}}
  \]
  for all $i,j$.
\end{lemma}

\begin{proof}
  We write $\phi(D_{\{1,2\}}) = D_{\{i,j\}}$ and $\phi(D_{\{1,3\}}) =
  D_{\{k,l\}}$. Assume that $\{i,j\} \cap \{k,l\} = \emptyset$. Then
  by \lref{1lMonomials} a monomial~$t$ of degree~$r$ occurs in an
  element of $I$ such that~$t$ is divisible by $D_{\{i,j\}}^2
  D_{\{k,l\}}$. By the hypothesis, $\phi^{-1}(t)$ also occurs in an
  element of $I$. But $\phi^{-1}(t)$ is divisible by $D_{\{1,2\}}^2
  D_{\{1,3\}}$, contradicting \lref{1lMonomials}. This argument shows
  that if the index sets of two $D_{\{\nu,\mu\}}$'s intersect, then
  the same is true for their images under~$\phi$. This will be used
  several times during the proof. Here, after possibly reordering the
  index sets (recall that we do not assume $i < j$ or $k < l$) we
  obtain $i = l$.  Thus $\phi(D_{\{1,3\}}) = D_{\{i,k\}}$. Now we
  write $\phi(D_{\{1,4\}}) = D_{\{m,p\}}$ and conclude, as above, that
  $\{m,p\} \cap \{i,j\} \ne \emptyset$ and $\{m,p\} \cap \{i,k\} \ne
  \emptyset$. Assume, by way of contradiction, that $i \notin
  \{m,p\}$. Then $\{m,p\} = \{j,k\}$, so $\phi(D_{\{1,4\}}) =
  D_{\{j,k\}}$. By \lref{1lMonomials} a monomial~$t$ of degree~$r$
  occurs in an element of $I$ such that~$t$ is divisible by
  $D_{\{i,j\}} D_{\{i,k\}} D_{\{j,k\}}$. Then $\phi^{-1}(t)$ also
  occurs in a polynomial from $I$, but $\phi^{-1}(t)$ is divisible by
  $D_{\{1,2\}} D_{\{1,3\}} D_{\{1,4\}}$. This contradicts
  \lref{1lMonomials}. Hence our assumption was false and we conclude
  that $i \in \{m,p\}$, so with suitable renumbering
  $\phi(D_{\{1,4\}}) = D_{\{i,m\}}$.
  
  Replacing~4 by any other index between~4 and~$n$, we conclude that
  $\phi(D_{\{1,\mu\}}) = D_{\{i,\pi(\mu)\}}$ with $\pi$ a permutation
  from $S_n$ (where we may assign $\pi(1) = i$). Now take $\nu,\mu \in
  \{2 \upto n\}$ with $\nu \ne \mu$. Writing $\phi(D_{\{\nu,\mu\}}) =
  D_{\{x,y\}}$, we conclude that $\{x,y\} \cap \{i,\pi(\mu)\} \ne
  \emptyset$ and $\{x,y\} \cap \{i,\pi(\nu)\} \ne \emptyset$. But
  assuming $i \in \{x,y\}$ (after renumbering $i = x$, say) leads to
  the contradiction $\phi(D_{\{\nu,\mu\}}) = D_{\{i,y\}} =
  \phi(D_{\{1,\pi^{-1}(y)\}})$. Hence $\{x,y\} =
  \{\pi(\nu),\pi(\mu)\}$ and therefore $\phi(D_{\{\nu,\mu\}}) =
  D_{\{\pi(\nu),\pi(\mu)\}}$, which concludes the proof.
\end{proof}

\subsection{Most $n$-point configurations are reconstructible from
  distances} \label{1sReconstructible}

In this section $K$ is a field of characteristic not equal to~2 (e.g.,
$K = \RR$ or $K = \CC$) and $V$ is an $m$-dimensional vector space
over $K$ equipped with a non-degenerate symmetric bilinear form
$\langle \cdot,\cdot \rangle$. Let $G = \Or(V) \subseteq \GL(V)$ be
the orthogonal group given by this form. The following proposition is
folklore.

\begin{prop} \label{1pGram}
  Let $v_1 \upto v_n$, $w_1 \upto w_n \in V$ be vectors with
  \[
  \langle v_i,v_j \rangle = \langle w_i,w_j \rangle \quad \text{for
    all} \quad i,j \in \{1 \upto n\}.
  \]
  Set $r := \min\{n,m\}$. If some $(r \times r)$-minor of the Gram
  matrix $\left(\langle v_i,v_j\rangle\right)_{i,j = 1 \upto n} \in
  K^{n \times n}$ is non-zero, then there exists a $g \in G$ such that
  $w_i = g(v_i)$ for all~$i$.
\end{prop}

\begin{proof}
  After renumbering we may assume that $A := \left(\langle
    v_i,v_j\rangle\right)_{i,j = 1 \upto r}$ is invertible. In
  particular, $v_1 \upto v_r$ are linearly independent. By the
  hypothesis, the same holds for $w_1 \upto w_r$, and $v_i \mapsto
  w_i$ gives an isomorphism between $\bigoplus_{i=1}^r K v_i$ and
  $\bigoplus_{i=1}^r K w_i$ which respects the form. By Witt's
  extension theorem there exists a $g \in G$ with $g(v_i) = w_i$ for
  $i \leq r$. This concludes the proof for $n \leq m$. Now assume $n >
  m$ and take an index $i > m$. There exist $\alpha_1 \upto \alpha_m
  \in K$ such that $v_i = \sum_{j=1}^m \alpha_j v_j$. So for $1 \leq k
  \leq m$ we have $\langle v_k,v_i \rangle = \sum_{j=1}^m \langle
  v_k,v_j \rangle \cdot \alpha_j$. It follows that
  \[
  \begin{pmatrix}
    \alpha_1 \\
    \vdots \\
    \alpha_m
  \end{pmatrix}
  = A^{-1}
  \begin{pmatrix}
    \langle v_1,v_i \rangle \\
    \vdots \\
    \langle v_m,v_i \rangle
  \end{pmatrix}.
  \]
  By the hypothesis, it follows that $w_i$ can be expressed as a
  linear combination of $w_1 \upto w_m$ {\em with the same
    coefficients\/}. Therefore
  \[
  w_i = \sum_{j=1}^m \alpha_j w_j = \sum_{j=1}^m \alpha_j g(v_j) =
  g(v_i).
  \]
\end{proof}

We come to the main theorem of this section. We assume that $K$, $V$,
and~$m$ are as above. We write $V^n$ for the direct sum of~$n$ copies
of $V$, so an $n$-point configuration is an element from
$V^n$. $K[V^n]$ is the ring of polynomials on $V^n$.

\begin{theorem} \label{1tReconstructible}
  Let $n$ be a positive integer with $n \leq 3$ or $n \geq m +2$. Then
  there exists a non-zero polynomial $f \in K[V^n]$ such that every
  $n$-point configuration $(P_1 \upto P_n)$ with $f(P_1 \upto P_n) \ne
  0$ is reconstructible from distances.
\end{theorem}

\begin{proof}
  The cases $n = 1$ or $m = 0$ are trivial. The case $m = 1$ will be
  proved in \sref{2sVolumes} (see \tref{2tReconstructible}). Therefore
  we may assume that $2 \leq n \leq 3$ or $2 \leq m \leq n - 2$.

  Take indeterminates $D_{\{i,j\}}$ indexed by sets $\{i,j\} \subset
  \{1 \upto n\}$ with $i \ne j$ and form the matrix
  \begin{equation} \label{1eqD2}
    {\cal D} := \left(D_{\{i,j\}} - D_{\{i,n\}} -
      D_{\{j,n\}}\right)_{i,j = 1 \upto n-1},
  \end{equation}
  where we set $D_{\{i,i\}} := 0$ as usual. If $2 \leq m \leq n-2$,
  let $I$ be the ideal of $(m+1) \times (m+1)$-minors of $\cal D$.
  Each permutation $\pi \in S_n$ induces a permutation $\phi_\pi$ of
  the $D_{\{i,j\}}$ by $\phi_\pi\left(D_{\{i,j\}}\right) =
  D_{\{\pi(i),\pi(j)\}}$. Let $H \leq S_{\binom{n}{2}}$ be the
  subgroup containing all the $\phi_\pi$, and let $\cal T$ be a set of
  left coset representatives of $H$, so we have a disjoint union
  \[
  S_{\binom{n}{2}} = \bigcup_{\psi \in {\cal T}}^. \psi H.
  \]
  We may assume that $\id \in {\cal T}$. \lref{1lPermutations} says
  that for every $\psi \in {\cal T} \setminus \{\id\}$ there exists an
  $F_\psi \in I$ such that $\psi(F_\psi) \notin I$. Set $F_1 :=
  \prod_{\psi \in {\cal T} \setminus \{\id\}} \psi(F_\psi)$. If, on
  the other hand, $2 \leq n \leq 3$, set $F_1 := 1$. In either case,
  set $r := \min\{n-1,m\}$ and let $F_2$ be a non-zero $(r \times
  r)$-minor of $\cal D$ (e.g., choose the first~$r$ rows and
  columns). Now set $F := F_1 F_2$.
  
  We choose a basis of $V \cong K^m$ such that $\langle \cdot,\cdot
  \rangle$ takes diagonal form, so $\langle (\xi_1 \upto
  \xi_m),(\eta_1 \upto \eta_m) \rangle = \sum_{k=1}^m a_k \xi_k
  \eta_k$ with $a_k \in K \setminus \{0\}$. Let $x_{i,j}$ be further
  indeterminates ($i = 1 \upto n$, $j = 1 \upto m$), so $K[V^n]$ can
  be identified with $K[x_{1,1} \upto x_{n,m}]$. Let
  $\map{\Phi}{K[D_{\{1,2\}} \upto D_{\{n-1,n\}}]}{K[x_{1,1} \upto
  x_{n,m}]}$ be the homomorphism of algebras given by $D_{\{i,j\}}
  \mapsto \sum_{k=1}^m a_k \left(x_{i,k} - x_{j,k}\right)^2$ (see
  \pref{1pRelations}(b)). Recall that $I$ is the kernel of $\Phi$.
  Since $\psi(F_\psi) \notin I$ for all $\phi \in {\cal T} \setminus
  \{\id\}$ and $F_2 \notin I$ (since each non-zero homogeneous element
  in $I$ has degree $> m$), we obtain that $f := \Phi(F) \ne 0$.

  Let $P_1 \upto P_n \in V$ such that $f(P_1 \upto P_n) \ne 0$, and
  let $d_{\{i,j\}} = \langle P_i - P_j,P_i - P_j\rangle$ be the
  distances. We have
  \begin{equation} \label{1eqNz}
    F\left(d_{\{1,2\}} \upto d_{\{n-1,n\}}\right) = f(P_1 \upto P_n)
    \ne 0.
  \end{equation}
  We wish to show that $P_1 \upto P_n$ form a reconstructible
  $n$-point configuration. Let $Q_1 \upto Q_n \in V$ be points with
  distances $d_{\{1,2\}}^\prime \upto d_{\{n-1,n\}}^\prime$ such that
  the distribution of distances coincides with that of the $P_i$. Then
  there exists a permutation~$\phi$ of the set ${\cal J} := \{\{i,j\}
  \subseteq \{1 \upto n\} \mid i \ne j\}$ (the index set of the D's)
  such that $d_{\{i,j\}}^\prime = d_{\phi(\{i,j\})}$. There exists a
  permutation $\pi \in S_n$ such that $\phi = \psi \circ \phi_\pi$
  with $\psi \in {\cal T}$. Thus
  \[
  d_{\psi(\{i,j\})} = d_{\{\pi^{-1}(i),\pi^{-1}(j)\}}^\prime
  \]
  for all $\{i,j\} \in {\cal J}$. Assume, by way of contradiction,
  that $\psi \ne \id$. Then $n \geq m+2$, since for $n \leq 3$ all
  permutations of $\cal J$ are induced from permutations from
  $S_n$. Clearly $\phi_{\pi^{-1}}$ preserves the ideal $I$, hence
  $F_\psi \in I$, implies $\phi_{\pi^{-1}}(F_\psi) \in I$. Therefore
  \[
  F_\psi\left(d_{\{\pi^{-1}(1),\pi^{-1}(2)\}}^\prime \upto
  d_{\{\pi^{-1}(n-1),\pi^{-1}(n)\}}^\prime\right) =
  \left(\phi_{\pi^{-1}}(F_\psi)\right)(d_{\{1,2\}}^\prime \upto
  d_{\{n-1,n\}}^\prime) = 0,
  \]
  and hence
  \[
  \left(\psi(F_\psi)\right)(d_{\{1,2\}} \upto d_{\{n-1,n\}}) =
  F_\psi\left(d_{\psi(\{1,2\})} \upto d_{\psi(\{n-1,n\})}\right) = 0,
  \]
  contradicting~\eqref{1eqNz}. It follows that $\psi = \id$, so
  $d_{\{i,j\}}^\prime = d_{\{\pi(i),\pi(j)\}}$ for all $i,j$. We have
  to show that there exists $g \in \AO(V)$ with $Q_i =
  g(P_{\pi(i)})$. For this purpose we may assume that $\pi$ is the
  identity. By applying a shift with a vector from $V$ we may further
  assume $P_n = Q_n = 0$. It follows from \eref{1eqSD} that the Gram
  matrices $\left(\langle P_i,P_j \rangle\right)_{i,j = 1 \upto n-1}$
  and $\left(\langle Q_i,Q_j \rangle\right)_{i,j = 1 \upto n-1}$
  coincide. Moreover,~\eqref{1eqNz} implies that an $(r \times
  r)$-minor of the Gram matrices is non-zero. Now \pref{1pGram} yields
  the desired result.
\end{proof}

\begin{rem} \label{1rInter}
  For $4 \leq n \leq m + 1$ (the range not covered by
  \tref{1tReconstructible}), no relations exist between the distances
  $d_{\{i,j\}}$ of an $n$-point configuration. If $K$ is algebraically
  closed, it follows from the surjectiveness of the categorical
  quotient (see \mycite[Theorem~3.5(ii)]{Newstead} or
  \mycite[Lemma~2.3.2]{Derksen:Kemper}) that for any given values for
  the $d_{\{i,j\}}$ there exists an $n$-point configuration which has
  these distances. Therefore in this case no $n$-point configuration
  is reconstructible from distances, with the possible exception of
  configurations where many of the distances are the same. It is not
  entirely clear whether the same holds for $K$ not algebraically
  closed (e.g. $K = \RR$), since in this case the categorical quotient
  is no longer surjective. As an example, for $K = \RR$ the distances
  must satisfy triangle inequalities. Nevertheless, we expect that
  also for $K = \RR$ and $4 \leq n \leq m + 1$, all $n$-point
  configurations lying in some dense open subset are not
  reconstructible from distances.
\end{rem}

%SECTION ADDED BY MIMI
\subsection{Symmetric $n$-point configurations} \label{1sSymmetric}

The reconstructibility test 
provided by Theorem \ref{1tReconstructible} 
fails for a variety of point configurations, 
including all those with repeated distances. 

\begin{lemma} \label{failure for repetitions}
 Let $P_1,\ldots, P_n\in V$ with $2\leq m \leq n-2$ 
 and consider $f$, 
 the polynomial function constructed in the proof of 
 Theorem \ref{1tReconstructible}. 
 If the pairwise distances between the $P_i$'s are not all distinct
 then $f(P_1,\ldots ,P_n)=0$.
\end{lemma}

\begin{proof}
 Denote by $d_{\{ i,j \}}$ the distance between $P_i$ and $P_j$.
 Assume that there exists $i_1,j_1,i_2, j_2$ 
 with $\{i_1,j_1 \}\neq \{ i_2, j_2\}$ 
 such that $d_{\{ i_1,j_1\}}=d_{\{i_2,j_2\}}$.
 Consider the permutation $\varphi\in S_{\binom{n}{2}}$ 
 which permutes $\{ i_1,j_1\}$ and $\{ i_2,j_2\}$
 and leaves all the other pairs $\{i,j \}$ unchanged.
 Observe that there does not exist $\pi\in S_n$ 
 such that $\varphi \{ i,j \}=\{\pi(i),\pi(j) \}$, for all $i,j=1,\ldots ,n$.
  %[Perhaps this is clear. But, just in case, here is the proof.
  %If there were such a $\pi$, we would have
  %$\{\pi(i_1),\pi(j_1) \} = \{ i_2,j_2\}$,
  %$\{\pi(i_2),\pi(j_2) \} = \{ i_1,j_1\}$,
  %and $\{\pi(i),\pi(j) \} = \{ i,j\}$ 
  %for all $(i,j)\neq (i_1,j_1),(i_2.j_2)$.
  %Consider $j_3$ such that $j_3 \neq i_1, j_1, j_2$ 
  %(this is possible since $n>3$).
  %We have $\{ \pi (i_1), \pi (j_3)\}= \{ i_1,j_3\}$ 
  %and $\{ \pi (i_1), \pi (j_2)\}= \{ i_1,j_2\}$. 
  %This implies that $\pi (i_1)=i_1$ 
  %since $i_1$ is the only index that appears in both equations. 
  %Therefore we also have that $\pi (j_2)=j_2$. 
  %But this contradicts the fact that 
  %$\{ \pi (i_1), \pi (j_1)\}= \{ i_2,j_2\}$.]
 Therefore, there exists $\psi \in {\cal T}\setminus \{ \id \}$
 and $\varphi_\pi\in H$ induced by a permutation $\pi \in S_n$
 such that $\varphi= \psi \circ \varphi_\pi$.

 Let $F_\psi$ be any polynomial with $F_\psi \in I$ 
 such that $\psi (F_\psi)\notin I$. 
   %\m{I think $\Psi$ from
   %\pref{1pRelations}(b) is not quite what we need here, so I changed
   %this a bit.}
% Let $\Psi$ be the homomorphism of algebra 
% defined in \pref{1pRelations}(b).
 We have $d_{\psi \{ i,j\}}= d_{\{ \pi^{-1}(i),\pi^{-1}(j) \}}$,
 for all $i,j=1,\ldots,n$.
 This means that
 \begin{eqnarray*}
 0 &=& F_\psi \left(d_{\{\pi^{-1}(1),\pi^{-1}(2)\}} \upto
   d_{\{\pi^{-1}(n-1),\pi^{-1}(n)\}}\right),
 \text{ since }F_\psi\in I,\\
   &=& F_\psi \left(d_{\psi(\{1,2\})} \upto
     d_{\psi(\{n-1,n\})}\right), \\
   &=& \psi F_{\psi}\left(d_{\{1,2\}} \upto d_{\{n-1,n\}}\right).
 \end{eqnarray*}
 So one of the factors of $f(P_1,\ldots, P_n)$ is zero 
 and the conclusion follows.
\end{proof}

\begin{cor}
 If an $n$-point configuration $P_1,\ldots,P_n$
 with $2\leq m \leq n-2$
 has a non-trivial symmetry,
 i.e.~if there exists $g\in \AO(V)$ and $\pi \in S_n\setminus \{\id\}$ 
%g could act trivially but s can not be the id,
%because if s is id, then g acts must act trivially, 
%so g and s don't change the point configuration
%and thus the symmetry is trivial.
 such that 
 \[
 (g\cdot P_1,\ldots,g\cdot P_n)=(P_{\pi(1)},\ldots,P_{\pi(n)}),
 \]
 then the polynomial function $f$ constructed in the proof 
 of Theorem \ref{1tReconstructible}
 is such that $f(P_1,\ldots,P_n)=0$.
\end{cor}

\begin{proof}
 By the previous lemma, it is sufficient to show that 
 there exists $\{i_1,j_1 \}\neq \{i_2,j_2 \}$
 such that $d_{\{ i_1,j_1\}}=d_{\{i_2,j_2\}}$.
 Since $\pi \neq \id$, 
 there exists $i_0$ such that $\pi(i_0)\neq i_0$.
 We have $g\cdot P_{i} =P_{\pi(i)}$, for all $i$'s, 
 so by invariance of the distance under $\AO(V)$, 
 this means that $d_{\{ i_0,j\}}=d_{\{\pi(i_0),\pi(j)\}}$ for all $j$'s.
 Therefore $i_1=i_0$, $i_2=\pi(i_0)$ $j_2=\pi(j_1)$
 and any  $j_1 \neq i_0,\pi(i_0)$ will do the trick.
\end{proof}

This does not mean that no symmetric $n$-point configuration is
reconstructible from distances.  Indeed a square is a counterexample
for $n=4$ (see \exref{1exRhombus} below).  We now give a reconstructibility test which does not
exclude all point configurations with repeated distances.
  %\m{Now I've
  %replaced your \pref{1pDoubleCosets} by another version. Reason: A
  %decomposition by cosets of $H \Sigma_{\underline{P}}$ as you suggest
  %doesn't exist, since $H \Sigma$ is not a subgroup. Please tell me
  %your opinion on the new version. Is it too technical?}

%Given a point configuration $P_1,\ldots ,P_n$ 
%whose pairwise distances $d_{\{ i,j \}}$ may not all be distinct,
%consider the set $\Sigma_{P_1,\ldots, P_n}$ 
%of all permutations  $\phi \in S_{\binom{n}{2}}$
%such that $d_{\phi\{i,j \} }=d_{\{ i,j \}}$, for all $i,j=1,\ldots,n$.

\begin{prop} \label{1pDoubleCosets}
  Let $P_1 \upto P_n \in V$ be points in an $m$-dimensional vector
  space ($2 \le m \le n-2$) over a field $K$ of characteristic not~2
  equipped with a non-degenerate symmetric bilinear form $\langle
  \cdot,\cdot \rangle$. Set $d_{\{i,j\}} := \langle P_i - P_j,P_i -
  P_j\rangle$, and assume that the matrix $\left(d_{\{i,j\}} -
    d_{\{i,n\}} - d_{\{j,n\}}\right)_{i,j = 1 \upto n-1}$ has rank~$m$
  (the ``generic'' rank). Let $G \le S_{\binom{n}{2}}$ be the subgroup
  of all permutations $\phi$ with $d_{\phi(\{i,j\})} = d_{\{i,j\}}$
  for all~$i$, $j$. (In fact, $G$ may be replaced by any smaller
  subgroup.) Moreover, let $H \le S_{\binom{n}{2}}$ be the subgroup of
  all $\phi_\pi$ with $\pi \in S_n$, given by $\phi_\pi(\{i,j\}) =
  \{\pi(i),\pi(j)\}$. Consider a set ${\cal T} \subset
  S_{\binom{n}{2}}$ of double coset representatives with respect to
  $G$ and $H$, e.i.,
  \[
  S_{\binom{n}{2}} = \bigcup_{\psi \in {\cal T}}^. G \psi H.
  \]
  Assume that $\id \in {\cal T}$, and for each $\psi \in {\cal T}
  \setminus \{\id\}$ choose $F_\psi \in I \setminus \psi^{-1}(I)$
  (where $I$ is the ideal occurring in Lemmas~\ref{1lMonomials}
  and~\ref{1lPermutations}), which is possible by
  \lref{1lPermutations}. If
  \[
  \left(\psi(F_\psi)\right)\left(d_{\{1,2\}} \upto
    d_{\{n-1,m\}}\right) \ne 0
  \]
  for all $\psi \in {\cal T} \setminus \{\id\}$, then $(P_1 \upto
  P_n)$ is reconstructible from distances.
\end{prop}

\begin{proof}
  Since the proof is almost identical to the one of
  \tref{1tReconstructible}, we will be very brief here to avoid
  repetitions. Let $Q_1 \upto Q_n \in V$ be points with (squared)
  distances $d_{\{i,j\}}^\prime$ such that $d_{\{i,j\}}^\prime =
  d_{\phi(\{i,j\})}$ with $\phi \in S_{\binom{n}{2}}$. Write $\phi =
  \rho \circ \psi \circ \phi_\pi$ with $\rho \in G$, $\psi \in {\cal
    T}$, and $\pi \in S_n$.  Then
  \[
  d_{\psi(\{i,j\})} = d_{(\rho \circ \psi)(\{i,j\})} = d_{(\phi \circ
    \phi_{\pi^{-1}})(\{i,j\})} = d_{\phi(\{\pi^{-1}(i),\pi^{-1}(j)\})} 
  = d_{\{\pi^{-1}(i),\pi^{-1}(j)\}}^\prime,
  \]
  where the first equality follows from the definition of $G$. As in
  the proof of \tref{1tReconstructible}, we conclude from this that
  $\psi = \id$, so $d_{\{i,j\}}^\prime = d_{(\rho \circ
    \phi_\pi)(\{i,j\})} = d_{\{\pi(i),\pi(j)\}}$ for all~$i$, $j$. The
  rest of the proof proceeds as for \tref{1tReconstructible}.
\end{proof}

\begin{ex} \label{1exRhombus}
  In this example we show that all rhombi are reconstructible from
  distances. Consider a rhombus in $K^2$ with sides of length~$a$ and
  diagonals of length~$b$ and~$c$ (see Figure~\ref{1fRhombus}), so
\Figure{
  \fbox{
    \begin{picture}(130,90)
      \put(5,45){\point}
      \put(-1,38){\attach{P_1}}
      \put(65,75){\point}
      \put(58,78){\attach{P_2}}
      \put(65,15){\point}
      \put(58,9){\attach{P_4}}
      \put(125,45){\point}
      \put(123,38){\attach{P_3}}
      \drawline(5,45)(65,75)
      \put(29,63){\attach{a}}
      \drawline(5,45)(65,15)
      \put(29,24){\attach{a}}
      \drawline(5,45)(125,45)
      \put(37,47){\attach{b}}
      \drawline(65,75)(65,15)
      \put(66,58){\attach{c}}
      \drawline(65,75)(125,45)
      \put(96,24){\attach{a}}
      \drawline(65,15)(125,45)
      \put(96,63){\attach{a}}
    \end{picture}
    }
  }
  {A rhombus}{1fRhombus}
  \[
  d_{\{1,2\}} = d_{\{2,3\}} = d_{\{3,4\}} = d_{\{1,4\}} = a,\
  d_{\{1,3\}} = b,\ \text{and} \ d_{\{2,4\}} = c.
  \]
  We assume that $a$, $b$, and $c$ are all non-zero. If we order the
  2-sets in $\{1 \upto 4\}$ as $\{1,2\}$, $\{1,3\}$, $\{1,4\}$,
  $\{2,3\}$, $\{2,4\}$, $\{3,4\}$, then the ``symmetry group'' $G$
  from \pref{1pDoubleCosets} is generated by the permutations $(1,3)$
  and $(1,3,4,6)$, and $G$ is isomorphic to $S_4$. The image $H$ of
  the embedding of $S_4$ into $S_6$ is generated by $(2,4)(3,5)$ and
  $(1,4,6,3)(2,5)$. It turns out that there are two double cosets in
  this case:
  \[
  S_6 = G H \stackrel{.}{\cup} G \psi H,
  \]
  where $\psi$ can be chosen as $\psi = (1,2)$. Since $m = 2$ and $n =
  4$, we have only one generating relation, which is the determinant
  of the matrix $\cal D$ defined in~\eqref{1eqD2}. Choose this
  determinant as the polynomial $F_\psi$. Assume that the rhombus is
  not reconstructible. By \pref{1pDoubleCosets} this implies
  $\left(\psi(F_\psi\right)\left(d_{\{1,2\}} \upto d_{\{3,4\}}\right)
  = 0$. We obtain
  \[
  a \left((a-b)^2 + c(c - b - 2 a)\right) = 0.
  \]
  We have $b + c = 4 a$. (This is Pythagoras' theorem, and it also
  follows from $F_\psi(d_{\{1,2\}} \upto d_{\{3,4\}}) = b c (b + c - 4
  a)$.) Substituting this into the above relation yields
  \[
  3 a (a - b) (c - a) = 0.
  \]
  Since $a \ne 0$, this implies $a = b$ or $a = c$ (here we need to
  assume that $\ch(K) \ne 3$), and by interchanging the roles of~$b$
  and~$c$ we may assume $a = b$. But this means that our rhombus has
  in fact a bigger symmetry group $\tilde{G}$ generated by the
  permutations $(1,2)$ and $(1,2,3,4,6)$. But now we see that $S_6 =
  \tilde{G} H$, so there is only the trivial double coset. It follows
  from \pref{1pDoubleCosets} that the rhombus is in fact
  reconstructible from distances.
  
  The computations for this example were done using the computer
  algebra system Magma~[\citenumber{magma}].
\end{ex}

\subsection{Locally reconstructible $n$-point configurations}

In this section,  $V$ is an $m$-dimensional vector space
over $K$ equipped with a non-degenerate symmetric bilinear form
$\langle \cdot,\cdot \rangle$.
We now concentrate 
on the local characterization of $n$-point configurations.
So we assume that $V^n$
is equipped with a norm $\| \cdot \|$
and that the field $K$ 
is either ${\mathbb R}$ or ${\mathbb C}$.
This first proposition
addresses the problem of local reconstructibility
for configurations of points whose mutual distances are all distinct.

\begin{prop}
 \label{neighborhood reconstructibility}
 Let $r=\min(n-1,m)$.
 Suppose that an $n$-point configuration $P_1,\ldots, P_n\in V$ 
 is such that
 its distances are all distinct
 and
 its Gram matrix (defined as in~\eqref{1eqD2}) has rank~$r$.
 Then there exists a neighborhood $N$ of $(P_1,\ldots,P_n)\in V^n$
 such that 
 any two $n$-point configurations in $N$
 are in the same orbit under the action of $\AO(V)$
 if and only if their distribution of distances is the same.
\end{prop}

\begin{proof}
 The distribution of distances is invariant under $\AO(V)$
 so one direction of the statement is trivial.
 To prove the other direction,
 observe that
 a minor is a determinant, which is a polynomial function, 
 and therefore continuous. 
 So there exists a neighborhood 
 $U$ of $(P_1,\ldots,P_n)\in V^n$
 such that
 the Gram matrix of any $(Q_1,\ldots,Q_n)\in U$
 has a non-zero $r$-by-$r$ minor.

 Let us assume the contrary, 
 so there exist two sequences of $n$-point configurations
 $\{ Q_1^k,\linebreak \ldots,Q_n^k \}_{k=1}^\infty$ and 
 $\{ R_1^k,\ldots,R_n^k \}_{k=1}^\infty$ in $U$,
 both converging to $P_1,\ldots,P_n$,
 and a sequence of permutations $\{ \varphi_k \}_{k=1}^\infty$, 
 such that for every $k$,
 $Q_1^k,\ldots,Q_n^k$ and $R_1^k,\ldots,R_n^k$ 
 are not in the same orbit under the action of $\AO(V)$
 but the distances 
 $d_{\{i,j \}}^{Q^k}=\langle Q_i^k-Q_j^k,Q_i^k-Q_j^k \rangle$
 are mapped to
 the distances
 $d_{\{i,j \}}^{R^k}=\langle R_i^k-R_j^k,R_i^k-R_j^k \rangle$
 by $\varphi_k$ 
 so   $d_{\{i,j \}}^{R^k}=d_{\varphi_k \{i,j \}}^{Q^k}$ 
 for all distinct $i,j=1,\ldots,n$.
 Since $S_{\binom{n}{2}}$ is finite, we may assume that
 $\varphi_k=\varphi$ is the same for every $k$.
 Taking the limit, we have
 \[
 \lim_{k\rightarrow \infty} d_{\{i,j \}}^{R^k} = 
 \lim_{k\rightarrow \infty} d_{\varphi\{i,j \}}^{Q^k}, 
 \text{ for all distinct } i,j=1,\ldots,n.
 \]
 By continuity of the distance,
 this implies that for any distinct $i,j=1,\ldots,n$
 the distance
 $d_{\{i,j \}}=\langle P_i-P_j,P_i-P_j \rangle$
 is equal to the distance
 $d_{\{\bar{i},\bar{j}\}}
 =\langle P_{\bar{i}}-P_{\bar{j}},P_{\bar{i}}-P_{\bar{j}}\rangle$
 where $\{\bar{i},\bar{j} \}=\varphi \{i,j \}$.
 Since all the $d_{\{i,j \}}$ are distinct, then $\varphi=\id$
 and thus  $d_{\{i,j \}}^{R^k}=d_{ \{i,j \}}^{Q^k}$ 
 for every distinct $i,j=1,\ldots,n$ and every $k$.
 By Proposition \ref{1pGram}, 
 this implies that
 $Q_1^k,\ldots,Q_n^k$
 and 
 $R_1^k,\ldots,R_n^k$
 are in the same orbit relative to $\AO(V)$, for every $k$
 which contradicts our hypothesis, and the conclusion follows.
\end{proof}

The following proposition
addresses the problem of local reconstructibility
for $n$-point configurations in general.

\begin{prop}
 \label{local reconstructibility}
 Let $r=\min(n-1,m)$.
 Suppose that an $n$-point configuration $P_1,\ldots, P_n\in V$ 
 is such that
 its Gram matrix (defined as in~\eqref{1eqD2}) has rank~$r$.
 Then there exists an $\epsilon >0$
 such that if the norm
 $ \| (Q_1,\ldots,Q_n)-(P_1,\ldots,P_n) \| < \epsilon $
 for some $n$-point configuration $Q_1,\ldots,Q_n \in V$
 with the same distribution of distances as that of $P_1,\ldots, P_n$,
 then $Q_1,\ldots,Q_n$ and  $P_1,\ldots, P_n$
 are in the same orbit relative to $\AO(V)$.
\end{prop}

\begin{proof}
 Again, by continuity, there exists a neighborhood $U$
 of $(P_1,\ldots,P_n)\in V^n$
 such that
 the Gram matrix of any $(Q_1,\ldots,Q_n)\in U$
 has a non-zero $r$-by-$r$ minor.
 Let us assume the contrary so
 there exists a sequence of $n$-point configurations 
 $\{ Q_1^k,\ldots,Q_n^k \}_{k=1}^\infty \subset U$
 converging to $P_1,\ldots,P_n$,
 and a sequence of permutations $\{ \varphi_k \}_{k=1}^\infty$, 
 such that none of the
 $Q_1^k,\ldots,Q_n^k$ 
 are in the same orbit as $P_1,\ldots,P_n$ under the action of $\AO(V)$
 but the distances 
 $d_{\{i,j \}}=\langle P_i-P_j,P_i-P_j \rangle$
 are mapped to
 the distances
 $d_{\{i,j \}}^{Q^k}=\langle Q_i^k-Q_j^k,Q_i^k-Q_j^k \rangle$
 by $\varphi_k$ 
 so   $d_{\varphi_k\{i,j \}}=d_{ \{i,j \}}^{Q^k}$ 
 for all $i,j=1,\ldots,n$ $i\neq j$.
 Again we may assume that
 $\varphi_k=\varphi$ is the same for every $k$.
 Taking the limit, we obtain that $ d_{\varphi\{i,j \}} = 
 \lim_{k\rightarrow \infty} d_{\{i,j \}}^{Q^k}$, 
 for all distinct  $i,j=1,\ldots,n$.
 By continuity of the distance,
 this implies that $ d_{\varphi\{i,j \}} = d_{\{i,j \}}$.
 Therefore, 
 $d_{ \{i,j \}}=d_{ \{i,j \}}^{Q^k}$
 for every $k$ and every distinct $i,j=1,\ldots,n$.
 By Proposition \ref{1pGram}, 
 this implies that
 $Q_1^k,\ldots,Q_n^k$ and $P_1,\ldots,P_n$
 are in the same orbit relative to $\AO(V)$ for every $k$,
 which contradicts our hypothesis, and the conclusion follows.
\end{proof}

When $V={\mathbb R}^m$, 
(the case that interests us the most for applications)
we can actually 
drop the requirement on the Gram matrix
based on the following refinement of Proposition \ref{1pGram}.

\begin{lemma}
 Let $G = \Or(V) \subseteq \GL(V)$ be
 the orthogonal group given by the form $\langle \cdot ,\cdot \rangle$.
 Let $v_1 \upto v_n$, $w_1 \upto w_n \in {\mathbb R}^m$ be vectors with
  \[
  \langle v_i,v_j \rangle = \langle w_i,w_j \rangle \quad \text{for
   all} \quad i,j \in \{1 \upto n\}.
  \]
 Then there exists a $g \in G$ such that
 $w_i = g(v_i)$ for all~$i$.
\end{lemma}

\begin{proof}
 Observe that since $V={\mathbb R}$, the rank of the Gram matrix 
 $\left(\langle v_i,v_j\rangle\right)_{i,j = 1 \upto n} $
 is equal to the dimension of the vector space spanned by $v_1,\ldots,v_n$.
 (This is {\em not} true over the complex field.)
 So we may assume, after relabeling, 
 that $v_1,\ldots,v_\rho$ with $\rho\geq 1$, are linearly independent.
 By hypothesis, the same is true for $w_1,\ldots,w_\rho$.
 By Proposition \ref{1pGram},
 there exists $g\in G$ such that $g (v_i)=w_i$, for all $i=1,\ldots,\rho$.

 For any $k$ such that $n\geq k>\rho$, 
 there exists $\alpha_1,\ldots,\alpha_\rho$ such that
 $v_k=\sum_{j=1}^\rho \alpha_j v_j$.
 So for $1 \leq k \leq \rho$ 
 we have 
 $\langle v_k,v_i\rangle=\sum_{j=1}^\rho\langle v_i,v_j\rangle\cdot\alpha_j$. 
 It follows that
  \[
  \begin{pmatrix}
    \alpha_1 \\
    \vdots \\
    \alpha_\rho
  \end{pmatrix}
  =\left( \left( \langle v_i,v_j\rangle\right)_{i,j=1\upto\rho }\right)^{-1}
  \begin{pmatrix}
    \langle v_1,v_i \rangle \\
    \vdots \\
    \langle v_\rho,v_i \rangle
  \end{pmatrix}.
  \]
  By the hypothesis, $w_i$ can be expressed as a
  linear combination of $w_1 \upto w_m$ 
  {\em with the same coefficients\/}. 
  Therefore
  \[
  w_i = \sum_{j=1}^m \alpha_j w_j = \sum_{j=1}^m \alpha_j g(v_j) =
  g(v_i).
  \]
\end{proof}

\begin{cor}
 For any $n$-point configuration $P_1,\ldots, P_n\in {\mathbb R}^m$
 whose distances are all distinct,
 there exists a neighborhood $N$ 
 of $(P_1,\ldots,P_n)\in \left( {\mathbb R}^m\right)^n$
 such that 
 any two $n$-point configurations in $N$
 are in the same orbit under the action of $\AO(V)$
 if and only if their distribution of distances is the same.
\end{cor}

\begin{cor}
 For any $n$-point configuration $P_1,\ldots, P_n\in {\mathbb R}^m$
 there exists an $\epsilon >0$ 
 such that if the norm
 $ \| (Q_1,\ldots,Q_n)-(P_1,\ldots,P_n) \| < \epsilon $
 for some $n$-point configuration $Q_1,\ldots,Q_n \in V$
 with the same distribution of distances as that of $P_1,\ldots, P_n$,
 then $Q_1,\ldots,Q_n$ and  $P_1,\ldots, P_n$
 are in the same orbit relative to $\AO(V)$.
\end{cor}

\section{Reconstruction from volumes} \label{2sVolumes}

Given $n$ points $P_1 \upto P_n \in \RR^2$ in a plane, we may consider
all areas $A_{i,j,k}$ of triangles spanned by three of these points
$P_i$, $P_j$, and $P_k$. Clearly these areas are preserved by the
action of all translations and all linear maps with determinant $\pm
1$. As in the preceding section, we can consider the {\em
distribution} of areas, and ask whether an $n$-point configuration is
reconstructible from this distribution up to the above action and
permutations of the points. Again we will generalize this to
configurations of points $P_i$ lying in $K^m$, with $K$ a field
and~$m$ any dimension. Since we are interested in invariants which are
preserved by all linear maps with determinant $\pm 1$, it makes sense
to consider volumes of $m$-simplices spanned by $m + 1$ points
$P_{i_0} \upto P_{i_m}$. These volumes are conveniently expressed by
the determinants
\begin{equation} \label{2eqVolumes}
  a_{i_0 \upto i_m} := \det\left(P_{i_1} - P_{i_0} \upto
  P_{i_m} - P_{i_0}\right)
\end{equation}
(where the $P_i$ are takes to be column vectors). The determinants are
really the ``signed volumes'', so we need to consider them up to
signs, which is equivalent to taking squares. This discussion leads to
the following definition.

\begin{defi} \label{2dReconstructible}
  Let $K$ be a field and $n > m$ positive integers. For an $n$-point
  configuration $P_1 \upto P_n \in K^m$ form the ``volumes'' $a_{i_0
  \upto i_m}$ as in \eqref{2eqVolumes} and the polynomial
  \[
  V_{P_1 \upto P_n}(X) = \prod_{1 \leq i_0 < \cdots < i_m \leq n}
  \left(X - a_{i_0 \upto i_m}^2\right).
  \]
  ($V_{P_1 \upto P_n}(X)$ encodes the distribution of volumes.) An
  $n$-point configuration $P_1 \upto P_n \in K^m$ is called
  \df{reconstructible from volumes} if the following holds: If $Q_1
  \upto Q_n$ is another $n$-point configuration with $V_{Q_1 \upto
  Q_n}(X) = V_{P_1 \upto P_n}(X)$, then there exist a permutation $\pi
  \in S_n$, a linear map $\phi \in \GL_m(K)$ with $\det(\phi) = \pm
  1$, and a vector $v \in K^m$ such that
  \[
  Q_i = \phi\left(P_{\pi(i)} + v\right)
  \]
  for all $i = 1 \upto n$.
\end{defi}

\begin{rem} \label{2rSmall}
  \begin{enumerate}
  \item If we are working in the plane, i.e., $m = 2$, we will of
    course speak of reconstructibility from {\em areas} instead of
    volumes.
  \item For $m = 1$, the above concept of reconstructibility from
  volumes coincides with reconstructibility from distances introduced
  in \dref{1dReconstructible}. \remend
  \end{enumerate}
  \renewcommand{\remend}{}
\end{rem}

\subsection{Non-reconstructible configurations} \label{2sCounter}

Again the first issue is to find configurations which are not
reconstructible from volumes. Our main interest will be
two-dimensional real space. A computation in
Magma~[\citenumber{magma}] yields that in $\RR^2$ all 4-point
configurations are reconstructible from volumes. For $n = 5$ we obtain
counterexamples (whose construction also involved Magma
computations). One of the simplest of these is given in
Figure~\ref{2fFive}.

\Figure{
  \fbox{
    \begin{picture}(130,90)
      \dottedline[\circle*{0.4}]{2}(0,9)(130,9)
      \dottedline[\circle*{0.4}]{2}(0,33)(130,33)
      \dottedline[\circle*{0.4}]{2}(0,57)(130,57)
      \dottedline[\circle*{0.4}]{2}(0,81)(130,81)
      \dottedline[\circle*{0.4}]{2}(5,89)(5,1)
      \dottedline[\circle*{0.4}]{2}(29,89)(29,1)
      \dottedline[\circle*{0.4}]{2}(53,89)(53,1)
      \dottedline[\circle*{0.4}]{2}(77,89)(77,1)
      \dottedline[\circle*{0.4}]{2}(101,89)(101,1)
      \dottedline[\circle*{0.4}]{2}(125,89)(125,1)
      \put(5, 33){\point}
      \put(6,27){\attach{P_1}}
      \put(29, 33){\point}
      \put(30,27){\attach{P_2}}
      \put(29, 57){\point}
      \put(20,59){\attach{P_3}}
      \put(77, 57){\point}
      \put(68,59){\attach{P_4}}
      \put(125, 57){\point}
      \put(116,59){\attach{P_5}}
    \end{picture}
  }
  \hspace{11mm}
  \fbox{
    \begin{picture}(130,90)
      \dottedline[\circle*{0.4}]{2}(0,21)(130,21)
      \dottedline[\circle*{0.4}]{2}(0,45)(130,45)
      \dottedline[\circle*{0.4}]{2}(0,69)(130,69)
      \dottedline[\circle*{0.4}]{2}(5,89)(5,1)
      \dottedline[\circle*{0.4}]{2}(29,89)(29,1)
      \dottedline[\circle*{0.4}]{2}(53,89)(53,1)
      \dottedline[\circle*{0.4}]{2}(77,89)(77,1)
      \dottedline[\circle*{0.4}]{2}(101,89)(101,1)
      \dottedline[\circle*{0.4}]{2}(125,89)(125,1)
      \put(29, 21){\point}
      \put(29,15){\attach{Q_1}}
      \put(53, 21){\point}
      \put(53,15){\attach{Q_2}}
      \put(53, 45){\point}
      \put(43,48){\attach{Q_3}}
      \put(53, 69){\point}
      \put(43,72){\attach{Q_4}}
      \put(101, 69){\point}
      \put(91,72){\attach{Q_5}}
    \end{picture}
  }
}
{Two 5-point configurations with the same distribution of areas}
{2fFive}

We put the points on a grid of length~1. The two configurations in
Figure~\ref{2fFive} lie in different orbits of $S_5 \times
\AGL_2(\RR)$, since in the first configuration all points lie on two
parallel lines, which is not the case in the second configuration. But
the signed areas $a_{i,j,k}$ are as follows:
\begin{center}
  \begin{tabular}{c||c|c|c|c|c|c|c|c|c|c}
    & $a_{1,2,3}$ & $a_{1,2,4}$ & $a_{1,2,5}$ & $a_{1,3,4}$ &
    $a_{1,3,5}$ & $a_{1,4,5}$ & $a_{2,3,4}$ & $a_{2,3,5}$ &
    $a_{2,4,5}$ & $a_{3,4,5}$ \\ \hline
    P & 1 & 1 & 1 & -2 & -4 & -2 & -2 & -4 & -2 & 0 \\ \hline
    Q & 1 & 2 & 2 & 1 & -1 & -4 & 0 & -2 & -4 & -2
  \end{tabular}
\end{center}
So the distributions of areas coincide.

For $n = 6$ we get an even simpler example which is given in
Figure~\ref{2fSix}.

\Figure{
  \fbox{
    \begin{picture}(130,90)
      \dottedline[\circle*{0.4}]{2}(0,65)(130,65)
      \dottedline[\circle*{0.4}]{2}(0,25)(130,25)
      \dottedline[\circle*{0.4}]{2}(5,10)(5,80)
      \dottedline[\circle*{0.4}]{2}(45,10)(45,80)
      \dottedline[\circle*{0.4}]{2}(85,10)(85,80)
      \dottedline[\circle*{0.4}]{2}(125,10)(125,80)
      \put(5,65){\point}
      \put(45,65){\point}
      \put(125,65){\point}
      \put(5,25){\point}
      \put(45,25){\point}
      \put(125,25){\point}
    \end{picture}
  }
  \hspace{11mm}
  \fbox{
    \begin{picture}(130,90)
      \dottedline[\circle*{0.4}]{2}(0,65)(130,65)
      \dottedline[\circle*{0.4}]{2}(0,25)(130,25)
      \dottedline[\circle*{0.4}]{2}(5,10)(5,80)
      \dottedline[\circle*{0.4}]{2}(45,10)(45,80)
      \dottedline[\circle*{0.4}]{2}(85,10)(85,80)
      \dottedline[\circle*{0.4}]{2}(125,10)(125,80)
      \put(5,65){\point}
      \put(45,65){\point}
      \put(125,65){\point}
      \put(5,25){\point}
      \put(85,25){\point}
      \put(125,25){\point}
    \end{picture}
  }
}
{Two 6-point configurations with the same distribution of areas}
{2fSix}

The configurations in Figure~\ref{2fSix} lie in different orbits of
$S_6 \times \AGL_2(\RR)$ since the first configuration has three
connecting vectors between points which are equal and the second one
has not. But it is easy to see that the configurations have the same
distribution of areas. Moreover, we can add an arbitrary number of
points on the upper dotted line in both configurations to obtain pairs
of $n$-point configurations with equal distributions of areas for $n
\ge 6$.

To get examples in dimension $m \ge 3$, one can embed the
two-dimensional examples given here into $m$-space and then add the $m
- 2$ points with coordinates $(0,0,1,0 \upto 0) \upto (0 \upto 0,1)$.

\subsection{Relation-preserving permutations} \label{2sRelation}

In this section $K$ is a field, $n$ and~$m$ are positive integers with
$n > m$, and $x_{i,j}$ are indeterminates ($1 \leq i \leq n$, $1 \leq
j \leq m$). For $1 \leq i_0 < \cdots < i_m \leq n$ we take further
indeterminates $A_{i_0 \upto i_m}$. Write $K[\underline{A}]$ for
the polynomial ring in the $A_{i_0 \upto i_m}$ and let $I
\subseteq K[\underline{A}]$ be the kernel of the map
\[
\mapl{\Phi}{K[\underline{A}]}{K[\underline{x}]}{A_{i_0 \upto
      i_m}}{\det\left(x_{i_j,k} - x_{i_0,k}\right)_{j,k = 1 \upto m}}.
\]
For $i_0 \upto i_m \in \{1 \upto n\}$ pairwise distinct, select the
permutation $\pi$ of the set $\{0 \upto m\}$ such that $i_{\pi(0)} <
i_{\pi(1)} < \cdots < i_{\pi(m)}$ and set
\begin{equation} \label{2eqAi}
  A_{i_0 \upto i_m} := \sgn(\pi) \cdot A_{i_{\pi(0)} \upto
  i_{\pi(m)}}.
\end{equation}

\begin{lemma} \label{2lRelations}
  \begin{enumerate}
  \item If $i_0 \upto i_{m+1} \in \{1 \upto n\}$ are pairwise
    distinct, then
    \[
    \sum_{k=0}^{m+1} (-1)^k A_{i_0 \upto i_{k-1},i_{k+1} \upto
    i_{m+1}} \in I.
    \]
  \item $I$ is generated by the polynomials $\sum_{k=0}^{m+1} (-1)^k
    A_{i_0 \upto i_{k-1},i_{k+1} \upto i_{m+1}}$ with $1 \leq i_0 <
    \cdots < i_{m+1} \leq n$ and by homogeneous polynomial of degree
    $> 1$ which only involve the $A_{n,i_1 \upto i_m}$ with $1 \leq
    i_1 < \cdots < i_m < n$.\footnote{The non-linear polynomials are
    the well-known Pl\"ucker relations, which we do not need to
    present here explicitly.}
  \item For $j \in \{1 \upto n\}$ the $A_{j,i_1 \upto i_m}$ with $1
  \leq i_1 < \cdots < i_m \leq n$, $i_k \ne j$, are linearly
  independent modulo~$I$.
  \end{enumerate}
\end{lemma}

\begin{proof}
  It is convenient to write $P_i$ for the (column) vector $(x_{i,1}
  \upto x_{i,m})^{\operatorname{T}}$, so for $i_0 \upto i_m \in
  \{1 \upto n\}$ in increasing order we have
  \begin{equation} \label{2eqA}
    \Phi\left(A_{i_0 \upto i_m}\right) = \det\left(P_{i_1} -
    P_{i_0} \upto P_{i_m} - P_{i_0}\right),
  \end{equation}
  which is equal to $\sum_{k=0}^m (-1)^k \det\left(P_{i_0} \upto
  P_{i_{k-1}},P_{i_{k+1}} \upto P_{i_m}\right)$. This shows
  that~\eqref{2eqA} is also valid if the $i_j$ are not increasing.
  \begin{enumerate}
  \item By~\eqref{2eqA} we have
    \begin{multline*}
      \Phi\left(A_{i_0 \upto i_m}\right) = \\
      \det\left(\strut\left(P_{i_1} - P_{i_{m+1}}\right) -
      \left(P_{i_0} - P_{i_{m+1}}\right) \upto \left(P_{i_m} -
      P_{i_{m+1}}\right) - \left(P_{i_0} - P_{i_{m+1}}\right)\right) =
      \\ \Phi\left(A_{i_{m+1},i_1 \upto i_m}\right) -
      \Phi\left(A_{i_{m+1},i_0,i_2 \upto i_m}\right) +- \cdots +
      (-1)^m \Phi\left(A_{i_{m+1},i_0 \upto i_{m-1}}\right) = \\
      \Phi\left(A_{i_0 \upto i_{m-1},i_{m+1}}\right) -+ \cdots +
      (-1)^m \Phi\left(A_{i_1 \upto i_m,i_{m+1}}\right).
    \end{multline*}
    This yields~(a).
  \item The relations between the $\Phi\left(A_{n,i_1 \upto
    i_m}\right)$ are known from classical invariant theory 
    (see \linebreak
    \mycite{Weyl} or \mycite{deConcini:Procesi}) to be the Pl\"ucker
    relations, which are homogeneous and non-linear. Let $J \subseteq
    K[\underline{A}]$ be the ideal generated by the linear relations
    given in~(b) and the Pl\"ucker relations. By~(a) we have $J
    \subseteq I$. Conversely, take $f \in I$. Using the linear
    relations from~(b), we can substitute every $A_{i_0 \upto i_m}$
    appearing in~$f$ by $\sum_{k=0}^m (-1)^k A_{n,i_0 \upto
    i_{k-1},i_{k+1} \upto i_m}$. In this way we obtain $g \in
    K[\underline{A}]$ with $f \equiv g \mod{J}$, and~$g$ only involves
    indeterminates $A_{i_0 \upto i_m}$ with $i_0 = n$. But $f \in I$
    implies $g \in I$, so~$g$ lies in the ideal generated by the
    Pl\"ucker relations. Thus $f \in J$.
  \item It follows from~(b) that the $\Phi\left(A_{n,i_1 \upto
    i_m}\right)$ with $1 \leq i_1 < \cdots < i_m < n$ are linearly
    independent. But the same argument can be made with any other
    index~$j$ instead of~$n$. This implies~(c). \endproof
  \end{enumerate}
  \renewcommand{\endproof}{}
\end{proof}

The next lemma shows that the linear relations given in
\lref{2lRelations} are the only ones of their kind.

\begin{lemma} \label{2lUnique}
  Let $l \in K[\underline{A}]$ be a non-zero linear combination of at
  most $m + 2$ of the indeterminates $A_{i_0 \upto i_m}$. Assume
  that all the coefficients in~$l$ are~1 or~-1, and $l \in I$. Then
  \begin{equation} \label{2eqLinrel}
    l = \sum_{k=0}^{m+1} (-1)^k A_{i_0 \upto i_{k-1},i_{k+1} \upto
    i_{m+1}}
  \end{equation}
  with $i_1 \upto i_{m+2} \in \{1 \upto n\}$ pairwise distinct.
\end{lemma}

\begin{proof}
  Take any $A_{i_0 \upto i_m}$ which occurs in~$l$. Define a
  homomorphism $\map{\phi}{K[\underline{A}]}{K[\underline{A}]}$ by
  sending each $A_{j_0 \upto j_m}$ with $i_0 \in \{j_0 \upto j_m\}$ to
  itself and by sending each $A_{j_0 \upto j_m}$ with $i_0 \notin
  \{j_0 \upto j_m\}$ to $\sum_{k=0}^m (-1)^k A_{i_0,j_0 \upto
  j_{k-1},j_{k+1} \upto j_m}$. \lref{2lRelations}(a) implies that
  $\phi(f) \equiv f \mod{I}$ holds for all $f \in
  K[\underline{A}]$. Thus $\phi(l) \in I$. But by
  \lref{2lRelations}(c) this implies $\phi(l) = 0$. But $A_{i_0 \upto
  i_m}$ occurs as a summand in $\phi(l)$ and must therefore be
  cancelled out by something. Hence a summand of the form $\pm
  A_{j_0,i_1 \upto i_m}$ with $j_0 \notin \{i_0 \upto i_m\}$ must
  occur in~$l$. The same argument can be applied to the other indices
  of $A_{i_0 \upto i_m}$, and we find summands $\pm A_{i_0 \upto
  i_{k-1},j_k,i_{k+1} \upto i_m}$ with $j_k \notin \{i_0 \upto i_m\}$
  in~$l$. We have already found $m + 2$ summands in~$l$, hence these
  are all summands.
  
  Now we apply the same argument to $A_{j_0,i_1 \upto i_m}$. Doing so
  we find that for each $k \in \{1 \upto m\}$ there must occur an
  indeterminate in~$l$ whose indices include all of $j_0,i_1 \upto
  i_{k-1},i_{k+1} \upto i_m$. Ruling out all other possibilities, we
  see that this indeterminate must be $A_{i_0 \upto i_{k-1},j_k,i_{k+1}
  \upto i_m}$, so $j_k = j_0$. Setting $i_{m+1} := j_0$, we find that
  up to the signs the summands of~$l$ are as claimed in the lemma.

  If $K$ has characteristic~2 then nothing has to be shown about signs
  and we are done. So assume $\ch(K) \ne 2$ and write $l' :=
  \sum_{k=0}^{m+1} (-1)^k A_{i_0 \upto i_{k-1},i_{k+1} \upto
  i_{m+1}}$. Assume that~$l$ is neither~$l'$ nor $-l'$. Since $l'$
  lies in $I$ by \lref{2lRelations}(a), the same is true for $(l +
  l')/2$. But $(l + l')/2$ is non-zero, has coefficients $\pm 1$, and
  has fewer than $m + 2$ summands. By the above discussion, this is
  impossible. Hence we conclude that $l = \pm l'$. Performing a
  permutation with sign $-1$ on the indices transforms $l'$ into
  $-l'$, so the case $l = -l'$ is also dealt with.
\end{proof}

The following proposition is analogous to \lref{1lPermutations}.
\newcommand{\eps}{\varepsilon}

\begin{prop} \label{2pPermutations}
  Let $\map{\phi}{K[\underline{A}]}{K[\underline{A}]}$ be an
  algebra-automorphism sending each $A_{i_0 \upto i_m}$ to $\pm
  A_{j_0 \upto j_m}$ for some $j_0 \upto j_m \in \{1 \upto
  n\}$ (where the signs may be chosen independently). If $\phi(I)
  \subseteq I$, then there exists $\pi \in S_n$ and $\eps \in \{\pm
  1\}$ such that for $1 \leq i_0 < \cdots < i_m \leq n$ we have
  \[
  \phi\left(A_{i_0 \upto i_m}\right) = \eps \cdot A_{\pi(i_0) \upto
  \pi(i_m)}.
  \]
\end{prop}

\begin{proof}
  If $n = m + 1$, there is only one indeterminate $A_{i_0 \upto i_m}$,
  so there is nothing to show. Hence we may assume that $n \geq m +
  2$. Set ${\cal M} := \left\{S \subset \{1 \upto n\} \mid |S| = m +
  1\right\}$. We have a bijection $\map{\psi}{\cal M}{\cal M}$ induced
  from $\phi$ by defining $\psi\left(\{i_0 \upto i_m\}\right) = \{j_0
  \upto j_m\}$ if $\phi\left(A_{i_0 \upto i_m}\right) = \pm A_{j_0
  \upto j_m}$. For $S = \{i_0 \upto i_m\} \in {\cal M}$ with $i_0 <
  \cdots < i_m$ we write $A_S := A_{i_0 \upto i_m}$, so $\phi(A_S) =
  \pm A_{\psi(S)}$. The bulk of the proof consists of constructing a
  permutation $\pi \in S_n$ such that
  \begin{equation}  \label{2eqPsiPi}
    \psi(S) = \pi(S)
  \end{equation}
  for all $S \in {\cal M}$, where the right-hand side means
  element-wise application of $\pi$.

  Take a subset $T \subseteq \{1 \upto n\}$ with $m + 2$ elements and
  write $T = \{i_0 \upto i_{m+1}\}$ with $i_0 < \cdots < i_{m+1}$. By
  \lref{2lRelations}(a) the polynomial $l = \sum_{k=0}^{m+1} (-1)^k
  A_{T \setminus \{i_k\}}$ lies in $I$, hence also $\phi(l) \in
  I$. But $\phi(l) = \sum_{k=0}^{m+1} \pm A_{\psi(T \setminus
  \{i_k\})}$. From \lref{2lUnique} we see that $\tilde{T} :=
  \bigcup_{k=0}^{m+1} \psi\left(T \setminus \{i_k\}\right)$ must have
  precisely $m + 2$ elements. Since each $\psi\left(T \setminus
  \{i_k\}\right)$ has $m + 1$ elements, there exists a map
  $\map{\pi_T}{T}{\tilde{T} \subseteq \{1 \upto n\}}$ with
  $\psi\left(T \setminus \{i_k\}\right) = \tilde{T} \setminus
  \{\pi_T(i_k)\}$. Since $\psi$ is injective this also holds for
  $\pi_T$, so $\pi_T(T) = \tilde{T}$. Thus for all $S \in {\cal M}$
  with $S \subset T$ we have
  \begin{equation} \label{2eqPsiPires}
    \psi(S) = \pi_T(S)
  \end{equation}
  (where the right-hand side means element-wise application of
  $\pi_T$).

  In the sequel we will make frequent use of the following rule: If
  two sets $S,S' \in {\cal M}$ have~$m$ elements in common, then also
  $\psi(S)$ and $\psi(S')$ share~$m$ elements. Indeed, there is a
  linear polynomial~$l$ of the type~\eqref{2eqLinrel} in which both
  $A_S$ and $A_{S'}$ occur. By \lref{2lRelations}(a), $l$ lies in $I$,
  hence also $\phi(l) \in I$. But $A_{\psi(S)}$ and $A_{\psi(S')}$
  occur in $\phi(l)$, hence $|\psi(S) \cap \psi(S')| = m$ by
  \lref{2lUnique}.

  Now take two subsets $T$,~$T' \subseteq \{1 \upto n\}$ with $|T| =
  |T'| = m + 2$ such that $S := T \cap T'$ has $m + 1$ elements. We
  will show that $\pi_T$ and $\pi_{T'}$ coincide on $S$. Write
  \[
  T = S \cup \{j\} \quad \text{and} \quad T' = S \cup \{k\}
  \]
  with $j,k \in \{1 \upto n\}$. For $l \in S$ set $S_l := T' \setminus
  \{l\}$, so $S_l \in {\cal M}$. Then $|S_l \cap \left(T \setminus
  \{l\}\right)| = m$ and $|S_l \cap S| = m$, so $\psi(S_l)$ shares~$m$
  elements with $\psi\left(T \setminus \{l\}\right) = \pi_T(T)
  \setminus \{\pi_T(l)\}$ and with $\psi(S) = \pi_T(S) = \pi_T(T)
  \setminus \{\pi_T(j)\}$. But $\psi(S_l)$ cannot be a subset of
  $\pi_T(T)$ since this would imply
  \[
  \psi(S_l) = \pi_T\left(\pi_T^{-1}\left(\psi(S_l)\right)\right) =
  \psi\left(\pi_T^{-1}\left(\psi(S_l)\right)\right),
  \]
  contradicting the injectiveness of~$\psi$, since $S_l \not\subseteq
  T$. It follows that $\psi(S_l) = \pi_T\left(T \setminus
  \{j,l\}\right) \cup \{r_l\}$ with $r_l \in \{1 \upto n\} \setminus
  \pi_T(T)$. We can write this slightly simpler as $\psi(S_l) =
  \pi_T\left(S \setminus \{l\}\right) \cup \{r_l\}$. On the other
  hand, we have $S_l \subset T'$, so
  \[
  \psi(S_l) = \pi_{T'}(S_l) = \pi_{T'}\left(S \setminus \{l\}\right)
  \cup \{\pi_{T'}(k)\}.
  \]
  Intersecting the resulting equality $\pi_T\left(S \setminus
  \{l\}\right) \cup \{r_l\} = \pi_{T'}\left(S \setminus \{l\}\right)
  \cup \{\pi_{T'}(k)\}$ over all $l \in S$ yields $\bigcap_{l \in S}
  \{r_l\} = \{\pi_{T'}(k)\}$. Thus $r_l = \pi_{T'}(k)$ independently
  of~$l$, and $\pi_T\left(S \setminus \{l\}\right) = \pi_{T'}\left(S
  \setminus \{l\}\right)$ for all $l \in S$. This shows that $\pi_T(l)
  = \pi_{T'}(l)$ for all $l \in S$, as claimed.

  We proceed by taking any two subsets $T$,~$T' \subseteq \{1 \upto
  n\}$ with $|T| = |T'| = m + 2$. We can move from $T$ to $T'$ by
  successively exchanging elements. Using the above result, we see
  that $\pi_T$ and $\pi_{T'}$ coincide on $T \cap T'$. Thus we can
  define $\map{\pi}{\{1 \upto n\}}{\{1 \upto n\}}$ such that for every
  subset $T \subseteq \{1 \upto n\}$ with $|T| = m + 2$ the
  restriction $\pi|_{_T}$ coincides with $\pi_T$. Now~\eqref{2eqPsiPi}
  follows from~\eqref{2eqPsiPires}, and it also follows that $\pi \in
  S_n$.

  Take $S \in {\cal M}$ and write $S = \{i_0 \upto i_m\}$ with
  $i_0 < \cdots < i_m$. The definition of $\psi$
  and~\eqref{2eqPsiPi} imply that
  \[
  \phi\left(A_{i_0 \upto i_m}\right) = \eps_S \cdot A_{\pi(i_0) \upto
  \pi(i_m)}
  \]
  with $\eps_S \in \{\pm 1\}$. We wish to show that $\eps_S$ does not
  depend on $S$. To this end, take $T \subseteq \{1 \upto n\}$ with
  $|T| = m + 2$ and write $T = \{i_0 \upto i_{m+1}\}$ with $i_0 <
  \cdots < i_{m+1}$. By \lref{2lRelations}(a), $l := \sum_{k=0}^{m+1}
  (-1)^k A_{i_0 \upto i_{k-1},i_{k+1} \upto i_{m+1}}$ lies in $I$,
  hence $\phi(l) \in I$. But
  \[
  \phi(l) = \sum_{k=0}^{m+1} (-1)^k \eps_{T \setminus \{i_k\}} \cdot
  A_{\pi(i_0) \upto \pi(i_{k-1}),\pi(i_{k+1}) \upto \pi(i_{m+1})}.
  \]
  \lref{2lUnique} implies that all $\eps_{T \setminus \{i_k\}}$
  coincide. This shows that if two sets $S$,~$S' \in {\cal M}$
  share~$m$ elements, then $\eps_S = \eps_{S'}$. But since we can move
  from any $S \in {\cal M}$ to any other $S' \in {\cal M}$ by
  successively exchanging elements, it follows that indeed all
  $\eps_S$ coincide. This completes the proof.
\end{proof}

\subsection{Most $n$-point configurations are reconstructible from
  volumes} \label{2sReconstructible}

In this section $K$ is a field and $V$ is an $m$-dimensional vector
space over $K$. The following proposition is well known.

\begin{prop} \label{2pSeparating}
  Let $v_1 \upto v_n$, $w_1 \upto w_n \in V$ be vectors with $n \geq
  m$, such that for all $1 \leq i_1 < \cdots < i_m \leq n$
  \[
  d_{i_1 \upto i_m} := \det\left(v_{i_1} \ldots v_{i_m}\right) =
  \det\left(w_{i_1} \ldots w_{i_m}\right).
  \]
  If at least one of the $d_{i_1 \upto i_m}$ is non-zero, then there
  exists a $\phi \in \SL(V)$ such that $w_i = \phi(v_i)$ for all~$i$.
\end{prop}

\begin{proof}
  After renumbering we may assume that $d_{1,2 \upto m}$ is
  non-zero. Hence $v_1 \upto v_m$ and\linebreak 
  $w_1 \upto w_m$ are linearly
  independent, and there exists a (unique) $\phi \in \SL(V)$ such that
  $w_i = \phi(v_i)$ for all $i \leq m$. Assume $n > m$ and take an
  index $i > m$. There exist $\alpha_1 \upto \alpha_m \in K$ such that
  $v_i = \sum_{j=1}^m \alpha_j v_j$. Indeed, by Cramer's rule we have
  $\alpha_j = (-1)^{n-j} d_{1 \upto j-1,j+1 \upto m,i}/d_{1 \upto
  m}$. By the hypothesis, it follows that $w_i$ can be expressed as a
  linear combination of $w_1 \upto w_m$ {\em with the same
    coefficients\/}. Therefore
  \[
  w_i = \sum_{j=1}^m \alpha_j w_j = \sum_{j=1}^m \alpha_j \phi(v_j) =
  \phi(v_i).
  \]
\end{proof}

We come to the main theorem of this section. We assume that $K$ is a
field, $V$ is an $m$-dimensional vector space over $K$, and $n > m$ is
an integer. We write $V^n$ for the direct sum of~$n$ copies of $V$, so
an $n$-point configuration is an element from $V^n$. $K[V^n]$ is the
ring of polynomials on $V^n$.

\begin{theorem} \label{2tReconstructible}
  There exists a non-zero polynomial $f \in K[V^n]$ such that every
  $n$-point configuration $(P_1 \upto P_n)$ with $f(P_1 \upto P_n) \ne
  0$ is reconstructible from volumes.
\end{theorem}

\begin{proof}
  Clearly we may assume $m > 0$. For indices $1 \leq i_0 < \cdots <
  i_m \leq n$, let $A_{i_0 \upto i_m}$ be an indeterminate, and for
  $i_0 \upto i_m \in \{1 \upto n\}$ pairwise distinct define $A_{i_0
  \upto i_m}$ as in~\eqref{2eqAi}. Let $I \subseteq K[\underline{A}]$
  be the kernel of the map $\map{\Phi}{K[\underline{A}]}{K[V^n]}$
  sending $A_{i_0 \upto i_m}$ to the polynomial $\Phi\left(A_{i_0
  \upto i_m}\right)$ with $\Phi\left(A_{i_0 \upto i_m}\right)\left(P_1
  \upto P_n\right) = \det\left(P_{i_1} - P_{i_0} \upto P_{i_m} -
  P_{i_0}\right)$ for $P_1 \upto P_n \in V$. Note that $I$ is
  precisely the ideal introduced at the beginning of
  \sref{2sRelation}.

  Let $G \subseteq \Aut_K\left(K[\underline{A}]\right)$ be the group
  of all automorphisms $\phi$ of $K[\underline{A}]$ sending each
  $A_{i_0 \upto i_m}$ to $\pm A_{j_0 \upto j_m}$ with $1 \leq j_0 <
  \cdots < j_m \leq n$.  For each permutation $\pi \in S_n$ and each
  $\eps \in \{\pm 1\}$ there is an automorphism $\phi_{\pi,\eps} \in
  G$ with $\phi_{\pi,\eps}\left(A_{i_0 \upto i_m}\right) = \eps \cdot
  A_{\pi(i_0) \upto \pi(i_m)}$. Let $H \leq G$ be the subgroup of all
  these $\phi_{\pi,\eps}$, and choose a set ${\cal T}$ of left coset
  representatives of $H$ in $G$ with $\id \in {\cal
  T}$. \pref{2pPermutations} says that for every $\psi \in {\cal T}
  \setminus \{\id\}$ there exists an $F_\psi \in I$ such that
  $\psi(F_\psi) \notin I$. Set $F := A_{n,1,2 \upto m} \cdot
  \prod_{\psi \in {\cal T} \setminus \{\id\}} \psi(F_\psi)$ and $f :=
  \Phi(F) \in K[V^n]$. $F \notin I$ implies that $f \ne 0$.
  
  Let $P_1 \upto P_n \in V$ such that $f(P_1 \upto P_n) \ne 0$, and
  for $1 \leq i_0 < \cdots < i_m \leq n$ let $a_{i_0 \upto i_m} =
  \det\left(P_{i_1} - P_{i_0} \upto P_{i_m} - P_{i_0}\right)$ be the
  ``signed volume''. We have
  \begin{equation} \label{2eqNz}
    F\left(\underline{a}\right) = f(P_1 \upto P_n) \ne 0.
  \end{equation}
  We wish to show that $P_1 \upto P_n$ form a reconstructible
  $n$-point configuration. Let $Q_1 \upto Q_n \in V$ be points and set
  $a_{i_0 \upto i_m}^\prime := \det\left(Q_{i_1} - Q_{i_0} \upto
  Q_{i_m} - Q_{i_0}\right)$. Assume that the distribution of volumes
  of $Q_1 \upto Q_n$ coincides with that of $P_1 \upto P_n$, i.e.,
  $V_{Q_1 \upto Q_n}(X) = V_{P_1 \upto P_n}(X)$. This means that up to
  signs the $a_{i_0 \upto i_m}^\prime$ are a permutation of the
  $a_{i_0 \upto i_m}$, so there exists a $\phi \in G$ such that for
  all $H \in K[\underline{A}]$ we have
  \begin{equation} \label{2eqH}
    \left(\phi(H)\right)\left(\underline{a}\right) =
    H\left(\underline{a}'\right).
  \end{equation}
  There exist $\pi \in S_n$ and $\eps \in \{\pm 1\}$ such that $\phi =
  \psi \circ \phi_{\pi,\eps}$ with $\psi \in {\cal T}$. By way of
  contradiction, assume that $\psi \ne \id$. Clearly
  $\phi_{\pi^{-1},\eps}$ preserves the ideal $I$, hence $F_\psi \in I$
  implies $H := \phi_{\pi^{-1},\eps}(F_\psi) \in I$. Therefore
  $H\left(\underline{a}'\right) = \left(\Phi(H)\right)\left(Q_1 \upto
  Q_n\right) = 0$, so~\eqref{2eqH} yields
  \[
  \left(\psi(F_\psi)\right)\left(\underline{a}\right) =
  \left(\phi(H)\right)\left(\underline{a}\right) =
  H\left(\underline{a}'\right) = 0,
  \]
  contradicting~\eqref{2eqNz}. It follows that $\psi = \id$, so $\phi
  = \phi_{\pi,\eps}$. We have to show that there exist $v \in V$ and
  $\psi \in \GL(V)$ with $\det(\psi) \in \{\pm 1\}$ such that $Q_i =
  \psi\left(P_{\pi(i)} + v\right)$ for all~$i$. For this purpose we
  may assume that $\pi$ is the identity. If $\eps = -1$, we apply an
  (arbitrary) linear map with determinant~-1 to $Q_1 \upto Q_n$. This
  will change all the signs of the $a_{i_0 \upto i_m}^\prime$. Hence
  we may assume that $\eps = 1$, so $\phi = \id$, and~\eqref{2eqH}
  implies $a_{i_0 \upto i_m}^\prime = a_{i_0 \upto i_m}$ for all index
  vectors $i_0 \upto i_m$. Since $a_{n,1,2 \upto m} \ne 0$ (this was
  the purpose of introducing $A_{n,1,2 \upto m}$ as a factor into
  $F$), \pref{2pSeparating} yields that there exists $\sigma \in
  \SL(V)$ such that $\sigma(P_i - P_n) = Q_i - Q_n$ for all $i \in \{1
  \upto n-1\}$. Setting $v := \sigma^{-1}(Q_n) - P_n$ gives the
  desired result $Q_i = \sigma(P_i + v)$ for $ i \in \{1 \upto n\}$.
\end{proof}

\begin{rem} \label{2rDoubleCosets}
  Everything that was said in \sref{1sSymmetric} about
  reconstructibility of configurations with symmetries carries over to
  reconstructibility from volumes. In particular, the analogue of
  \pref{1pDoubleCosets} holds.  Similarly, the analogues of
  Propositions~\ref{neighborhood reconstructibility} and \ref{local
    reconstructibility} concerning local reconstructibility are also
  true.
\end{rem}

\subsection{Combining distances and volumes} \label{2sCombine}

Taking another look at Figure~\ref{1fCounter}, one notices that
although the two configuration have the same distribution of
distances, their distributions of areas are different. This brings up
the idea to try to distinguish $n$-point configurations (up to the
action of $S_n \times \AO_m(K)$) by considering the distribution of
distances {\em and} the distribution of volumes. Could it be that by
combining these data we might be able to separate all orbits? The
following example shows that once again this is not the case. We take
the following 4-point configurations in $\RR^2$ (see
Figure~\ref{2fCombination}):
\[
\begin{array}{llll}
  P_1 = (0,0), & P_2 = (0,6), & P_3 = (6 \sqrt{2},0), & P_4 = (2
  \sqrt{2},-1), \\
  Q_1 = (0,0), & Q_2 = (0,6), & Q_3 = (6 \sqrt{2},0), & Q_4 = (2
  \sqrt{2},5).
\end{array}
\]
\Figure{
  \fbox{
    \begin{picture}(130,90)
      \put(18.33095245999999999999999999, 17.5){\point}
      \put(8,13){\attach{P_1}}
      \put(18.33095245999999999999999999, 83.5){\point}
      \put(8,79){\attach{P_2}}
      \put(111.6690475519999999999999999, 17.5){\point}
      \put(113,13){\attach{P_3}}
      \put(49.44365082399999999999999999, 6.5){\point}
      \put(51,0){\attach{P_4}}
      \drawline(18.33095245999999999999999999, 17.5)%
      (18.33095245999999999999999999, 83.5)
      \drawline(18.33095245999999999999999999, 17.5)%
      (111.6690475519999999999999999, 17.5)
      \drawline(18.33095245999999999999999999, 17.5)%
      (49.44365082399999999999999999, 6.5)
      \drawline(18.33095245999999999999999999, 83.5)%
      (111.6690475519999999999999999, 17.5)
      \drawline(18.33095245999999999999999999, 83.5)%
      (49.44365082399999999999999999, 6.5)
      \drawline(111.6690475519999999999999999, 17.5)%
      (49.44365082399999999999999999, 6.5)
    \end{picture}
  }
  \hspace{11mm}
  \fbox{
    \begin{picture}(130,90)
      \put(18.33095245999999999999999999, 17.5){\point}
      \put(7,13){\attach{Q_1}}
      \put(18.33095245999999999999999999, 83.5){\point}
      \put(7,79){\attach{Q_2}}
      \put(111.6690475519999999999999999, 17.5){\point}
      \put(113,13){\attach{Q_3}}
      \put(49.44365082399999999999999999, 72.5){\point}
      \put(53,73){\attach{Q_4}}
      \drawline(18.33095245999999999999999999, 17.5)%
      (18.33095245999999999999999999, 83.5)
      \drawline(18.33095245999999999999999999, 17.5)%
      (111.6690475519999999999999999, 17.5)
      \drawline(18.33095245999999999999999999, 17.5)%
      (49.44365082399999999999999999, 72.5)
      \drawline(18.33095245999999999999999999, 83.5)%
      (111.6690475519999999999999999, 17.5)
      \drawline(18.33095245999999999999999999, 83.5)%
      (49.44365082399999999999999999, 72.5)
      \drawline(111.6690475519999999999999999, 17.5)%
      (49.44365082399999999999999999, 72.5)
    \end{picture}
  }
}
{Two 4-point configurations with the same distribution of distances
  and the same distribution of areas}
{2fCombination}

It is easy to see that the two configurations lie in different orbits
of $S_4 \times \AO_2(\RR)$ (although they lie in the same orbit of
$S_4 \times \AGL_2(\RR)$). We obtain the following distances
$\sqrt{d_{i,j}}$ and signed areas $a_{i,j,k}$:
%\begin{center}
%  \begin{tabular}{c||c|c|c|c|c|c||c|c|c|c}
%    & $d_{1,2}$ & $d_{1,3}$ & $d_{1,4}$ & $d_{2,3}$ & $d_{2,4}$ &
%    $d_{3,4}$ & $a_{1,2,3}^2$ & $a_{1,2,4}^2$ & $a_{1,3,4}^2$ &
%    $a_{2,3,4}^2$ \\ \hline
%    P & 36 & 72 & 9 & 108 & 57 & 33 & 2592 & 288 & 72 & 1800\\ \hline
%    Q & 36 & 72 & 33 & 108 & 9 & 57 & 2592 & 288 & 1800 & 72
%  \end{tabular}
%\end{center}
\begin{center}
  \begin{tabular}{c||c|c|c|c|c|c||c|c|c|c}
    & $\sqrt{d_{1,2}}$ & $\sqrt{d_{1,3}}$ & $\sqrt{d_{1,4}}$ &
    $\sqrt{d_{2,3}}$ & $\sqrt{d_{2,4}}$ & $\sqrt{d_{3,4}}$ &
    $a_{1,2,3}$ & $a_{1,2,4}$ & $a_{1,3,4}$ & $a_{2,3,4}$ \\ \hline
    P & 6 & 6 $\sqrt{2}$ & 3 & 6 $\sqrt{3}$ & $\sqrt{57}$ &
    $\sqrt{33}$ & $-36 \sqrt{2}$ & $-12 \sqrt{2}$ & $-6 \sqrt{2}$ &
    $-30 \sqrt{2}$ \\ \hline
    Q & 6 & 6 $\sqrt{2}$ & $\sqrt{33}$ & 6 $\sqrt{3}$ & 3 &
    $\sqrt{57}$ &  $-36 \sqrt{2}$ & $-12 \sqrt{2}$ & $30 \sqrt{2}$ &
    $6 \sqrt{2}$
  \end{tabular}
\end{center}

\section*{Acknowledgments}
We thank Serkan Hosten and Greg Reid 
for inviting us to the Symbolic Computational
Algebra conference held in London, Ontario in 2002. This is where we
first met and started this project.

The idea of using distributions of invariants
in order to separate the orbits was
inspired by discussions 
of Mireille Boutin with David Cooper and Senem Velipasalar
regarding their work on indexation~[\citenumber{TasdizenVelipasalarCooper}].
This author is grateful to
the SHAPE lab of Brown University
for providing the environment for these discussions
and thus the motivation for this paper.

\addcontentsline{toc}{section}{References}

\bibliographystyle{mybibstyle} \bibliography{bib}

\bigskip

\begin{center}
\begin{tabular}{lll}
  Mireille Boutin & &Gregor Kemper \\
  Max Planck Institute  & & Technische Universit\"at M\"unchen \\
  Inselstra\ss e 22 & & Zentrum Mathematik - M11 \\
  D-04103, Leipzig & & Boltzmannstr. 3 \\
  Germany & & 85\,748 Garching \\
  & & Germany \\
  {\tt boutin$@$mis.mpg.de} & & {\tt kemper$@$ma.tum.de}
\end{tabular}
\end{center}

\end{document}